\documentclass[11 pt, reqno]{amsart}
\usepackage[left= 24mm, right= 24mm]{geometry}

\usepackage[english]{babel}

\usepackage{amsmath}
\usepackage{amssymb}
\usepackage{amsfonts}
\usepackage{mathtools}
\usepackage{mathrsfs}
\usepackage{dsfont} 
\usepackage{stmaryrd}
\usepackage{xcolor}

\usepackage{graphicx}
\usepackage{tikz}
\usepackage{tikz-cd}
\usetikzlibrary{arrows.meta, positioning, calc}

\usepackage{verbatim} 
\usepackage{bm}

\usepackage{hyperref}

\hypersetup{
pdftitle={pdfTitle},
pdfauthor={AA},
colorlinks=true,
urlcolor=red
}

\newcommand{\N}{\mathbb{N}}

\newcommand{\R}{\mathbb{R}}

\newcommand{\al}{\alpha}

\newcommand{\ga}{\gamma}
\newcommand{\eps}{\varepsilon}

\newcommand{\E}{\mathbb{E}}

\renewcommand{\l}{ \left( }
\renewcommand{\r}{ \right) }

\renewcommand{\lq}{\left[}
\renewcommand{\rq}{\right]}

\renewcommand{\lg}{\left\{}
\newcommand{\rg}{\right\}}

\usepackage[backend=biber,
style=alphabetic,]{biblatex}
\def \Xt  { \mathcal{X}(t) }
\def \Vt  { \mathcal{V}(t) }
\def \X { \mathcal{X} }
\def \V { \mathcal{V} }

\def \sgM { \mathcal{T} } 

\def \sgL { \mathcal{Q} } 

\renewcommand{\S}{\mathcal{S}}

\def \bfn {\eta}

\def \Inv { L }
\def \Int { \mathcal{I} }
\def \J {\mathcal{J}}

\DeclareMathAlphabet{\mathpzc}{OT1}{pzc}{m}{it}

\newtheorem{te}{Theorem}[section]
\newtheorem{defin}[te]{Definition}
\newtheorem{os}[te]{Remark}
\newtheorem{prop}[te]{Proposition}
\newtheorem{lem}[te]{Lemma}
\newtheorem{ex}[te]{Example}
\newtheorem{corollary}[te]{Corollary}

\newtheorem*{assumption*}{\assumptionnumber}
\providecommand{\assumptionnumber}{}
\makeatletter
\newenvironment{assumption}[2]
{%
	\renewcommand{\assumptionnumber}{(\textbf{#1#2})}%
	\begin{assumption*}%
		\protected@edef\@currentlabel{(\textbf{#1#2})}%
	}
	{%
	\end{assumption*}
}

\usepackage{enumitem}

\date{\today}

\begin{document}

\title[Non-Markovian random flights]{From Non-Markovian Random Flights and Non-local Linear Boltzmann Equations to Anomalous Diffusion}



\author[]{Lorenzo Facciaroni*}

\keywords{Non-Markovian dynamics, random flights, anomalous diffusion, non-local operators, linear Boltzmann equations, scaling limits \\
* Email: \texttt{lorenzo.facciaroni@uniroma1.it},
 Department of Statistical Sciences, Sapienza - University of Rome}
	\date{\today}

\begin{abstract}
    This paper deals with two models of random flights in $\R^d\times \S^{d-1}$ exhibiting non-Markovian features. Both of them are governed by non-local counterparts of the linear Boltzmann equation; emphasis is placed on their hydrodynamic limits.
    The first model falls in the class of semi-Markov random flights; in \cite{ricciuti2023semi} the authors prove that if the distribution of the flight times lies in the domain of attraction of a stable law, the hydrodynamic limit yields a superdiffusion with ballistic behaviour. A unique characterization is given of such flight times; this result is of independent interest in  theory of semi-Markov processes. We prove that the law of such limiting anomalous diffusion solves $\int_{\mathcal{S}^{d-1}} \l \l \partial_t + cv \cdot \nabla_x \r^\al \rho(x,t) - \frac{t^{-\al}}{\Gamma(1-\al)}\, \delta_0(x - cvt) \r \mu(dv) = 0$.
    The second model yields, instead, a subdiffusion; 
    Hilbert expansion techniques are used to recover the well-known fractional heat equation $\partial_t^\alpha \rho(x,t) = \Delta \rho(x,t)$ and its kinetic derivation as low-density limit of a Lorentz-type model is proved.
\end{abstract}

\maketitle
\tableofcontents

\section{Introduction and preliminaries}

\subsection{Context and vision}
Random flights generalize the concept of a discrete random walk into continuous directional space, initially introduced by Pearson \cite{pearson1905problem}. The most classical setting involves a particle, starting in a (possibly) random position $x\in\R^d$ with (possibly) random direction velocity $v\in\S^{d-1}$ and constant speed $c$, that moves along a straight line in $\mathbb{R}^d$; at random times, governed by a counting process, a new direction is instantly chosen according to a random distribution with support $\S^{d-1}$ \cite{de2012flying,de2020random,kolesnik2005planar,orsingher2007random,garra2014random-KG,stadje1989exact,watanabe1970convergence}. Under the hypothesis of i.i.d.~exponential flight times of parameter $\lambda$ and uniform directions, the resulting position-velocity process $ \l X,V\r = \lg (X(t),V(t)),\ t\geq0 \rg$ is said to be the isotropic Markovian random flight in $\R^d$. A relevant fact is that its joint probability density function $f(x,v,t)$ solves the linear Boltzmann equation (or Kolmogorov backward equation)
\begin{align}
    \partial_t f(x,v,t) = cv\cdot \nabla_x f(x,v,t) + \lambda \int_{\S^{d-1}} \l f(x,v',t) - f(x,v,t) \r \mu(dv'), \label{Introduzione: linear boltzmann equation}
\end{align}
where $\mu$ denotes the uniform density on $\S^{d-1}$. Here $\lambda > 0$ is the rate of the Poisson process and is precisely the mean number of direction changes per unit time. 
The case $d = 1$ has received considerable attention in the literature; here the particle can only move left or right and the law of the position process satisfies the telegraph equation 
\begin{align}
    \partial_t^2 p(x,t) + \lambda  \partial_t p(x,t) = c^2 \partial_x^2 p(x,t). \label{Introduzione equazione del telegrafo}
\end{align}
Under the so-called hydrodynamic limit ($c\to \infty,\ \lambda \to \infty$ under the constraint $c^2/\lambda \to 1$) the random position $X(t)$ converges in law to the Brownian motion, whose distribution is the fundamental solution of the heat equation
\begin{align}
    \partial_t \rho(x,t) = \frac{1}{2}\Delta_x \rho(x,t), \quad \rho(x,0) = \delta(x). \label{Fractional heat equation - watanabe ed altri}
\end{align}
Different settings for random flights have been considered in literature. The interested reader can consult, e.g., \cite{orsingher2016motion,pinsky1976isotropic} for random motion in non-Euclidean spaces, \cite[Section 5]{facciaroni2025random} for a model with random speed and \cite{de2012random} for random flights not necessarily isotropic. 

Random flights appear in the context of kinetic theory as a scaling limits of Lorentz models (in $d\in\{2,3\}$), i.e.~a model for the motion of a particle among randomly located scatterers, also known as obstacles. It was originally used to describe the transport of electrons through a conductor (see \cite{lorentz1904motion}). The non-Markovian character of this process makes it analytically intractable. 

In a pioneering contribution, Gallavotti \cite{Gallavotti1972RigorousTO} gave the first rigorous treatment of a Lorentz process: a particle colliding elastically with hard spheres of radius $R$, whose centers are distributed according to a Poisson point process of intensity $\rho>0$. He proved that in the Boltzmann-Grad limit ($\rho\to\infty$, $R\to0$, with $\rho R^{d-1}=\lambda$ fixed) the distribution of the particle converges to the solution of the linear Boltzmann equation. Spohn \cite{spohn1978lorentz} extended this result to more general distributions of obstacles, proving the convergence of the process. Intuitively, recollisions become negligible in the limit. The results of Gallavotti and Spohn, combined with the diffusive limit of the linear Boltzmann equation, can thus be summarized in the following scheme for a Lorentz gas with randomly distributed obstacles. Other settings, such as a periodic distribution of obstacles, inelastic collisions, or systems in which particles also collide with one another, are not covered.


\begin{center}
\begin{tikzpicture}[
    box/.style={rectangle, draw, rounded corners=6pt, thick,
                align=flush center, text width=2.8cm,
                minimum height=1.9cm, inner sep=5pt,
                label={[font=\fontsize{6pt}{7pt}\selectfont\scshape, text=black!55,
                        anchor=north, yshift=-5pt]north:#1}},
    arrow/.style={-{Latex[length=3mm]}, thick},
    lim/.style={font=\footnotesize, align=center, inner sep=2pt},
    node distance=2.5cm,
]

\node[box=Microscopic model, text width=3.2cm, font=\small] (micro)
  {\rule{0pt}{0.4cm}{\bfseries Newton's\\equations}};

\node[box=Mesoscopic model, text width=3.2cm, right=of micro, font=\small] (meso)
  {\rule{0pt}{0.4cm}{\bfseries Linear Boltzmann\\equation}};

\node[box=Macroscopic model, text width=3.2cm, right=of meso, font=\small] (macro)
  {\rule{0pt}{0.4cm}{\bfseries Diffusion\\equation}};

\draw[arrow] (micro) -- node[lim, above=6pt] {Boltzmann--Grad\\limit} (meso);
\draw[arrow] (meso)  -- node[lim, above=6pt] {Diffusive limit} (macro);
\end{tikzpicture}
\end{center}


In the physical literature, there are several context where random flights are used to describe
dynamics. An example is the description of motile bacteria, such as the run-and-tumble motion of Escherichia coli \cite{angelani2013averaged,bechinger2016active,berg2004coli}.
Additionally, random flight models find extensive applications in polymer physics to characterize macromolecular conformations \cite{degennes1979scaling,doi1986theory,chandrasekhar1943stochastic}, as well as in radiative transfer and neutron transport theory \cite{davison1957neutron}.

Standard diffusion models recalled above fail when applied to complex physical or biological processes characterized by non-local dynamics and memory effects. Empirical evidence frequently reveals anomalous signatures, such as a non-linear growth of the mean squared displacement, long-range correlations, heavy-tailed marginal distributions, and trapping effects \cite{scher1975anomalous}. To properly describe these phenomena, generalized frameworks that account for non-locality and non-Markovian behaviour are strictly required.

Recent contributions tried to handle this problem, introducing random flights whose distribution is governed by non-local equations \cite{angelani2024anomalous,ricciuti2023semi}. In the case of \cite{angelani2024anomalous}, the authors replaced the first order derivative in the linear Boltzmann equation \eqref{Introduzione: linear boltzmann equation} with the Caputo fractional derivative $\partial_t^\alpha,\ \al\in(0,1)$ (see Remark \ref{Osservazione: subordinatore stabile e Caputo fractional derivative} for the definition of the operator), thus introducing non-locality in time. Instead, in \cite[Section 5]{ricciuti2023semi} the authors introduced a random flight governed by the following fractional equation, for $\al\in(0,1)$,
\begin{align*}
    \l \partial_t - cv\cdot \nabla_x\r^{\al} q(x,v,t) - \frac{t^{-\al}}{\Gamma(1-\al)}h(x+vt, v) = \lambda \int_{\S^2} \l q(x,v',t) - q(x,v,t)\r \mu(dv').
\end{align*}
The latter random flight displays Mittag-Leffler flight times and Mittag-Leffler space-displacement. As a consequence, the free flight time
is a proper random variable  with infinite mean. The authors also show the convergence of the random position to a superdiffusion with mean square displacement of order $t^2$ for large $t$.

Regarding the Lorentz models, non-Markovian scaling limits have been discussed in the last years. In some of these models,  the gas is immersed in an external force field; the latter generates non-Markovianity in the joint vector position-velocity after the kinetic limit. The interest reader can consult \cite{bobylev2001liouville,desvillettes2007boltzmann,nota2024diffusion,nota2022two}. 
On the other hand, in \cite[Section 3]{facciaroni2025random} the authors show that non-Markovianity may arise just from the random distribution of obstacles, when this ceases to be Poisson. In the latter paper, the limiting random flight is fully non-Markovian, exhibiting a superdiffusive behaviour for large $t$ (the mean squared displacement goes like $t^{2-\al},\ \al\in(0,1)$), and has infinite mean free path.


Up to our knowledge, the first rigorous anomalous diffusion limit from a kinetic equation appeared (on ArXiv) in \cite{mellet2011fractional} and independently in \cite{jara2009limit} (in dimension $d = 1$). In both cases, the solution of the Boltzmann equation converges to the solution of a fractional heat equation $\partial_t \rho(x,t) = -(-\Delta_x)^\al\rho(x,t)$, hence involving the fractional Laplacian as non-local operator in the space variable. It is worth noting that in both cases the scaling is such that the time is accelerated (\cite[(2.21)]{jara2009limit} and \cite[(10)]{mellet2011fractional}).
Other non-local counterparts of heat equation appeared in literature recently, approximating models of Lorentz type (see \cite[Section 3, Section 5]{facciaroni2025random} and \cite{matthies2024fractional}).

\subsection{Summary of the main results}
An important fact to keep in mind is the following: both non-local operators $\partial_t^{\al}$ and $\l \partial_t - cv\cdot \nabla_x\r^{\al}$ are related, in a way that we shall discuss later, to the inverse of a $\al$-stable subordinator. We recall that a subordinator $H = \{H(t),\ t\geq0\}$ is a strictly increasing Lévy process; to each subordinator is associated a so-called Bernstein function $\phi$ such that $\E e^{-sH(t)} = e^{-t\phi(s)}$ and vice-versa.

Our primary goal is to  discuss non-Markovian random flights whose distribution is governed by non-local equations,
exploiting the theory of Bernstein functions \cite{schilling2012bernstein}; fractional operators of \cite{angelani2024anomalous} and \cite[Section 5]{ricciuti2023semi} are re-obtained for $\phi(x) = x^{\al}, \ \al\in(0,1)$.

In Section \ref{Ricciuti paragrafo} we discuss a generalization of the semi-Markov random flight presented in \cite[Section 5]{ricciuti2023semi} in detail. Here the governing (backward) equation becomes
\begin{align}
    \phi\l \partial_t - cv\cdot \nabla_x\r q(x,v,t) - \bar \nu (t)h(x+vt, v)  = \int_{\S^2} \l q(x,v',t) - q(x,v,t)\r \mu(dv'), \label{introduzione: ricciuti governing equation}
\end{align}
where the operator $\phi\l \partial_t - cv\cdot \nabla_x\r$ is defined in Section \ref{subsection: Ricciuti Toaldo governing equation}. 
Equation \eqref{introduzione: ricciuti governing equation} is proven to govern the random motion described as follows. A particle starts from a random position $x\in\R^{d}$ and is shot towards an arbitrary direction $v\in\S^{d-1}$. In the Markovian case $(X,V)$, the time required before the first change of velocity is $W\sim\mathrm{Exp}(\lambda)$; here, instead, we consider the random flight $\l X_1^\phi,V_1^\phi\r := \lg \l X_1^\phi(t),V_1^\phi(t)\r,\ t\geq0 \rg$ where $V_1^\phi(t) := V(L(t))$ and $X_1^\phi(t) := X(0)+ \int_0^tV_1^\phi(s)ds$, where $L$ is the inverse process of a subordinator $H$ with Bernstein function $\phi$. As a consequence, the distribution of each flight time (the random time before the first collision) ceases to be exponential, and follows a distribution of the form $H(W)$ with positive support.

We discuss the hyperbolic scaling of such random flight in the case of $\phi(x) = x^\al$; in Theorem \ref{prop limiting fractional equation} we show that the distribution $\rho(x,t)$ of the limiting superdiffusion solves the following fractional equation 
\begin{align*}
    \int_{\mathcal{S}^2} \l \l \partial_t + cv \cdot \nabla_x \r^\al \rho(x,t) - \frac{t^{-\al}}{\Gamma(1-\al)}\, \delta_0(x -cvt) \r \mu(dv) = 0,
\end{align*}
where $\mu$ is the uniform distribution on $\mathcal{S}^2$.
Such equation has already appeared in \cite[Section 5]{magdziarz2015limit}; the authors found it as limit equation of a Lévy walk with waiting times in a domain of attraction of an $\al$-stable law. However, they use a completely different approach, relying on the theory developed by Meerschaert and co-authors \cite{meerschaert2008triangular}. 

Finally, in \ref{Teorema: H(W) coincide con la classe delle distribuzioni geometricamente infinitamente divisibili} we prove the following characterization: the class of distributions $H(W)$ coincides with the class of geometrically infinitely divisible distribution; this result is the continuous-time counterpart of the one in \cite{pachon2018discrete}. Among the possible flight time distributions, we show in Theorem \ref{Teorema: dominio di attrazione stabile J = H(W)} that $H(W)$ lies in the domain of attraction of an $\al$-stable law iff $\bar \nu(t) $, the tail of the Lévy measure of $H$, is regularly varying of index $-\al$ as $t \to \infty$; for such distributions, \cite[Theorem 5.9]{ricciuti2023semi} applies and $\{X_1(t),\ t\geq0\}$ converges to a superdiffusion with mean square displacement of order $t^2$. These characterization \ref{Teorema: H(W) coincide con la classe delle distribuzioni geometricamente infinitamente divisibili} and \ref{Teorema: dominio di attrazione stabile J = H(W)} turn out to be of independent interest in probability theory, especially concerning theory of semi-Markov processes.

The model in \cite{angelani2024anomalous} is extended by introducing the time non-local kinetic equation 
\begin{align}
    \phi\l \partial_t \r f(x,v,t) = -cv\cdot \nabla_x f(x,v,t) + \lambda \int_{\mathcal{S}^{d-1}} \l f(x,v',t) - f(x,v,t) \r \mu(dv'), \label{Introduzione governing equation de gregorio}
\end{align}
where $\phi\l \partial_t \r$ is the generalized Caputo derivative (properly defined in  \eqref{operatore per semi markov}), acting on suitable functions.
Equation \eqref{Introduzione governing equation de gregorio} is the governing equation of the random flight $\l X_2^\phi,V_2^\phi\r  = \lg \l X_2^\phi(t),V_2^\phi(t) \r,\ t\geq0 \rg$ where $X_2^\phi$ is constructed as time-change of the position of the Markovian random flight with the inverse of the subordinator $H$ with Bernstein function $\phi$ (Section \ref{De Gregorio: governing equation equazione governante e proprietà distributizionali}).
We find the mean squared displacement, showing that it models anomalous diffusion phenomena.  
In Section \ref{A Hilbert expansion method formally leads to the time-fractional heat equation} we use the Hilbert expansion approach to find the limiting time-fractional diffusion equation 
\begin{align*}
     \partial_t^\al \rho(x,t) = \frac{1}{2} \Delta \rho(x,t).
\end{align*}
This method is complementary to the standard time-change techniques. 

An important result of this paper is the following: in Section \ref{paragrafo: De Gregorio: limite lorentz} we show that the non-Markovian random flight $\l X_1^\phi,V_1^\phi\r$ is an approximation, under a low-density limit, of a mechanical model of Lorentz type. 
Contrary to the current literature on Lorentz gas, we here consider a direction of research which can be characterized as follows: non-Markovian limits, and hence non-local equations, may arise from a time-change argument. 
This approach allows to model rigorously systems where the random distribution of obstacles does not suffice to explain the anomalous behaviour of particles (e.g.~trapping phenomena). During this trapping state, physical time progresses while a particle may remain in the same microscopic position.
This provides a rigorous microscopic-macroscopic bridge to justify non-local kinetic equations, which are often introduced only phenomenologically.

\subsection{Preliminaries}
In the present paper we shall use results from very different theories; hence, for the reader convenience, we summarize the main facts in few pages.
The unit sphere shall be denoted with $\S^{d-1} := \{ v\in \mathbb{R}^d: |v|_e^2=1\}$, where $|\cdot|_e$ denotes the Euclidean norm. 

\subsubsection{Random flights and Lorentz models}
Let a particle start at a point $x\in \mathbb{R}^d$ and move along a unit vector $v_0$ with random speed $\mathcal{C}_0$, for a random waiting time $J_1$.
After the time interval $J_1$, the particle undergoes scattering and continues its motion along a new unit vector $v_1$  with a new speed $\mathcal{C}_1$, for another random waiting time $J_2$. This process repeats indefinitely. In principle, the speeds $\mathcal{C}_i$ may be assumed either  positive constants or  positive random variables. We call $\S \subseteq \S^{d-1}$ the state space of the directions $\{v_i,\ i\in \N\}$. 

Let $\tau_n := \sum_{i = 1}^n J_i$ be the instant of the $n$-th change of direction for $n\in \N_+$, with the convention $\tau_0 : = 0$,  and let $\mathcal{N}=\{\mathcal{N}(t),\ t\in \N\}$ be the process which counts the number of changes up to time $t$, i.e. $\mathcal{N}(t) := \max \{n \geq 0: \tau _n \leq t\}$.
At time $t\geq0$, the particle has unit velocity vector $V(t)$ and is located at the position $X(t)$.
To summarize the above description, we now give the rigorous definition of the joint process $\lg (X(t), V(t)), t\geq0 \rg$.

\begin{defin} \label{Definizione transport process random flight generica}
    Let $\l \Omega, \mathcal{A}, \mathds{P}^{(x,v)} \r,\ v\in \S,\ x\in\R^d$, be a family of probability spaces.
    Let $ \{J_n\}_{n \geq 1} $, $\{v_i,\ i\in \N\}$, $\{\mathcal{C}_i,\ i\in \N\}$ and $\mathcal{N}$ be random variables as defined above.
    The joint process $ \lg \l X(t),  V(t)\r ,\ t\geq0\rg$, such that
    \begin{align} \label{Isotropic: definizione velocità transp process}
        V(t) &:= v_{\mathcal{N}(t)} \\   
    X(t)  &:= X(0) + \int_0^t \mathcal{C}_sV(s)  ds 
     = X(0) +  \sum _{i=1}^{\mathcal{N}(t)} \mathcal{C}_iJ_i v_{i-1}+ \mathcal{C}_{\mathcal{N}(t)}\ (t-\tau _{\mathcal{N}(t)})\  v_{\mathcal{N}(t)}, \label{Isotropic: definizione posizione transp process}
    \end{align}
     where $\mathds{P}^{x,v}(X(0)= x ,v_0 = v)=1$, is said to be a transport process or random flight.
\end{defin}
Such definition is very general. Directions $\{v_i,\ i\in \N\}$ may be dependent or not and may take values in a general state space; the choice of i.i.d.~directions with uniform law on a symmetric domain would lead to an isotropic random flight. Furthermore, one might assume dependence between directions $\{v_i,\ i\in \N\}$, waiting times $\lg J_n,\ n\in \N_+ \rg$ and speeds $\{\mathcal{C}_i,\ i\in \N\}$. Initial position and velocity may be taken random as well. Similar constructions of Definition \ref{Isotropic: definizione velocità transp process} may hold in non-Euclidean spaces (see e.g.~\cite{pinsky1976isotropic} for the isotropic case over Riemannian manifolds or \cite{orsingher2016motion} for the motion of a particle in the hyperbolic space); the latter models will not be discussed here. Scenarios involving non-constant speeds (e.g., particles under the influence of gravity or external force fields) fall outside this framework but are of interest.

Pay attention to the fact that if the speeds $\mathcal{C}_i$ are random, it would be different to consider the random flight $\l X(t), V(t)\r$ where $V(t) \in \S^{d-1}$ is a unit vector and the random flight $\l X(t), U(t)\r$ where $U(t) := \mathcal{C}_tV(t) \in \R^d$. The following assumption shall be used in some parts of the paper.
\begin{assumption}{A1}{} \label{Assumption on random flights}
    Assume $\{v_i,\ i\in \N\}$ i.i.d.~unit vectors uniformly distributed on $\S = \S^{d-1}$, independent of the waiting times $\lg J_n,\ n\in \N_+ \rg$ and constant speeds, i.e.~$\mathcal{C}_i = c \in \R_+$ almost surely for each $i\in \N$.
\end{assumption}
In the setting \ref{Assumption on random flights}, the random flight of Definition \ref{Definizione transport process random flight generica} is isotropic and is uniquely identified by the counting process $\mathcal{N}$. Here, the position increases linearly in time.
We shall use the notation 
\begin{align}
    \mathcal{J} h(x,v) := \int _{\S^{d-1}} h(x,v')\mu (dv'), \qquad h\in C_0 \l \R^d \times \S^{d-1} \r, \label{definizione operatore J di scattering}
\end{align}
for the scattering operator $\J$, where $\mu$ denotes the uniform distribution on $\S^{d-1}$. 
\begin{ex}

Let the waiting times $\{J_n, n\geq 1\}$ be i.i.d.~exponential random variables with mean $\lambda ^{-1}$. Then, under the Assumption \ref{Assumption on random flights} the resulting transport process is the Markovian isotropic transport process, as studied, for example, in \cite{monin1956statistical,orsingher2007random,stadje1989exact,watanabe1970convergence}. In this setting, the counting process $\mathcal{N}$ is a Poisson process with intensity $\lambda$ (the mean number of direction changes per unit time) and $\tau_n$ follows a Gamma distribution of parameters $n,\lambda$. 
In this case, we shall denote the random flight process of Definition \ref{Definizione transport process random flight generica} by 
\begin{align*}
   \l X^{(c,\lambda)}, V^{(c,\lambda)}\r := \lg \l X^{(c, \lambda)}(t), V^{(c, \lambda)}(t) \r,\ t\geq 0 \rg ,
\end{align*}
where the superscript $(c,\lambda)$ indicates that the process depends on the two parameters $c$ and $\lambda$. Let us consider initial data $(x,v)\in \R^d\times \S^{d-1}$.
The mean square displacement of the Markovian transport process (see e.g.~\cite[Equation (2.24)]{angelani2024anomalous}) is
    \begin{align}
        M (t) := \E \left| X^{(c,\lambda)}(t) - x \right|^2 = \frac{2c^2}{\lambda^2} \l \lambda t - 1 + e^{-\lambda t} \r. \label{Preliminaries: MSD del markovian random flight}
    \end{align}
Moreover, consider the family of operators $\lg \sgM_t^ {(c, \lambda)} \rg_{t\geq0}$ defined by
\begin{align}
    \sgM_t^ {(c, \lambda)}h (x,v) := \mathds{E}^{(x,v)}\lq h \l X^ {(c, \lambda)}(t), V^{(c, \lambda)}(t) \r \rq \qquad h\in C_0\l \R^d \times \S^{d-1}\r  \label{semigruppo Boltzmann}\end{align}
where $\mathds{E}^{(x,v)}$ denotes the expectation under $\mathds{P}^{x,v}$.
For $h$ sufficiently regular, one has that $t \mapsto \sgM_t^ {(c, \lambda)}h$ is the unique solution
to the linear (abstract) Boltzmann equation on $C_0 (\mathbb{R}^{d} \times \S^{d-1})$
    \begin{align} 
        \partial_t u(t)=  \mathcal{A} u(t), \quad u(0)= h,\label{Boltzmann equation}
\end{align}
where $\mathcal{A} :=c\, v\cdot \nabla  +\lambda (\mathcal{J}-\mathcal{I})$ and $\mathcal{I}$ is the identity operator.
This result is standard and the proof can be obtained using \cite[Theorem 1.9.2]{kolokoltsov2011markov}.
\end{ex}

\subsubsection{Subordinators and inverses} \label{preliminaries: subordinators and inverses}
A subordinator is a strictly increasing Lévy process, meaning with stationary and independent increments; recall that the subordinator does not have fixed discontinuities. 
The class of subordinator is in correspondence 1-1 with a special class of functions, known as Bernstein functions. A function $\phi : (0, \infty) \longrightarrow \R$ is said to be a Bernstein function if $\phi \in C^\infty\l (0, \infty) \r$, \ $\phi(\bfn) > 0$ for each $ \bfn > 0$ and $ (-1)^{n-1} \phi^{(n)}\l \bfn \r \geq0 \text{ for all } n \in \mathbb{N}_+,\ \bfn > 0$. Each Bernstein function has a unique integral representation,  known as Lévy-Khintchine representation.
Precisely, a function $\phi(\cdot)$ is a Bernstein function iff there exist constants $a,b \in\R_+$ and a measure $\nu$ on $(0, \infty)$ satisfying
$
\int_0^\infty (1 \wedge s) \, \nu(ds) < \infty$
such that
\begin{align}
\phi(\bfn) = a + b \bfn + \int_0^\infty \l1 - e^{-\bfn s}\r\, \nu(ds). \label{Rappresentazione Lévy Kintchine Bernstein function}
\end{align}
The link with subordinator is the following: each Bernstein function is the Laplace exponent of a subordinator, i.e.~to each Bernstein function $\phi(\cdot)$ there corresponds a subordinator $H$ such that $\mathbb{E}e^{-\eta H(t)}= e^{-t\phi(\eta)}$, and vice versa. From now on, we restrict to the case of $a = 0$, i.e.~the subordinator is not killed (does not explode in finite time); furthermore, if not specified, we shall also assume that $b=0$, i.e.~the subordinator does not have drift.


A common way to construct non-Markovian process is to use the inverse process of $H$. If $Y=\{Y(t),\ t\geq0\}$ a Markov process, letting $\Inv := \{\Inv(t),\, t \in\R_+\}$ denote the inverse of $H$, i.e.
    \begin{align*}
        \Inv(t) := \inf\lg s\in\R_+ \mid H(s) > t\rg, \qquad t\in\R_+,
    \end{align*}
one can define the (possibly) non-Markovian process $\lg Y\l \Inv(t) \r,\ t\in\R_+\rg$.
Given that $H$ is strictly increasing, $L$ is continuous.

As a relevant example, if $\mathcal{N} := \{\mathcal{N}(t),\ t\in\R_+\}$ is a Poisson process of intensity $\lambda > 0$,
then $\mathcal{N}^\phi := \{\mathcal{N}^\phi(t),\ t\in\R_+\}$ is a renewal counting process with i.i.d.~waiting times, where $\mathcal{N}^\phi(t) := \mathcal{N}\l L(t)\r,\ t\in\R_+$. Each waiting time $J$ of $\mathcal{N}^\phi$ has survival function of the form
\begin{align} \label{funzione di sopravvivenza}
    S_\lambda(t) = P(J > t) = P\l \mathcal{N}(L(t)) = 0 \r =\mathbb{E} e^{-\lambda L(t)},
\end{align}
as a consequence of a conditioning argument on $L(t)$.
Moreover, $S_\lambda(\cdot)$ is the eigenfunction of an integro-differential operator, known as generalized Caputo derivative, which, for a sufficiently regular function $g(\cdot)$,   has the form
\begin{align}
    \phi\l \partial_t \r  g(t) := \partial_t \int _0^t (g(\tau)-g(0) )  \overline\nu (t-\tau)d\tau, \label{operatore per semi markov}
\end{align}
where $\nu$ is the Lévy measure of $H$.
Then it is true that
    \begin{align} \label{equazione autovalori}
     \phi\l \partial_t \r   S_\lambda(t) = -\lambda S_\lambda(t), \qquad S_\lambda(0)=1.
    \end{align}
A standard reference for this theory is \cite{meerschaert2019stochastic}.
\begin{os} \label{Osservazione: subordinatore stabile e Caputo fractional derivative}
    If $H$ is the $\al$-stable subordinator, then $\phi(x) = x^{\al}$ and the operator \eqref{operatore per semi markov} becomes the more famous Caputo fractional derivative
    \begin{align}
        \partial_t^\al g(t) := 
        \begin{cases}
            \frac{1}{\Gamma(1-\al)} \partial_t \int_0^t \l g(\tau) - g(0) \r \l t-\tau\r^{-\al} d\tau & \al \in (0,1) \\
            g'(t) & \al = 1
        \end{cases}. \label{Caputo fractional derivative}
    \end{align}
    In this case, the renewal counting process $\{\mathcal{N}(L(t)),\ t\geq0\}$ has Mittag-Leffler waiting times and is known as the time-fractional Poisson process (see \cite{meerschaert2011fractional,ScalasGorenfloMainardi2004} for details). An important feature of the stable subordinator is that its Bernstein function is homogeneous of index $\al$, i.e.~$\phi(cx) = c^\al\phi(x),\ c \in \R$. For this reason,     applying the definition of Caputo derivative \eqref{Caputo fractional derivative} to a function $g(t) = h(t/\ga)$ leads, after a change of variable $u = \tau/\gamma$, to
    \begin{align*}
        \partial_t^\al h(t/\ga) = \partial_t \l \ga^{1-\al} \int_0^{t/\ga} \l g(u) - g(0)\r \frac{(t-\ga u)^{-\al}}{\Gamma(1-\al)}du \r = \ga^{-\al} \l \partial_t^\al h \r (t/\ga), \quad \ga\in\R_+.
    \end{align*}
\end{os}

\section{Non-Markovian transport processes: a semi-Markov approach}
\label{Ricciuti paragrafo}

In this section we are interested in a class of random flights showing non-Markovian features. We shall use the notation of Section \ref{preliminaries: subordinators and inverses}, so $H$ is a strictly increasing subordinator, with Bernstein function $\phi$ and inverse process $L$.

Starting from the general definition of random flights \ref{Definizione transport process random flight generica} under the hypothesis \ref{Assumption on random flights}, we are interested in the case where the counting process is a renewal counting process $\mathcal{N}^\phi$ as recalled in Section \ref{preliminaries: subordinators and inverses}. The resulting random flight appears as follows: a particle starts the motion at $x\in \R^d$ and moves along the unit vector $v_1$ with constant speed $c$, until it has a collision after a random waiting time $J_1$. Then the particle moves (again with speed $c$) along the unit direction $v_2$ for a random time $J_2$, and so on. 
We shall use notation $\l X_1^{\phi},V_1^{\phi} \r := \lg \l X_1^\phi(t), V_1^{\phi}(t)\r, \ t\geq0 \rg$, where 
\begin{align}
    V_1^{\phi}(t) := V^{(c,\lambda)}(L(t)), \qquad X_1^\phi(t) := X^{(c,\lambda)}(0) + c\int_0^{t} V_1^\phi (s) ds, \qquad t\geq0. \label{ricciuti toaldo definizione random flight esplicita}
\end{align}
As in Section \ref{paragrafo de gregorio section},  we do not make explicit the dependence on $c$ and $\lambda$. 
It is worth noting that here the position $X_1^{\phi}$ is still a linear function of real time $t$; when there are trappings, i.e.~$L(t)$ is experiencing a plateau, the particle moves in a straight line rather than remaining in the same position, accumulating displacement at that frozen velocity over real time.

\subsection{Governing equation and its hyperbolic scaling} \label{subsection: Ricciuti Toaldo governing equation}

In this section we discuss the governing equation for the law of the semi-Markov random flight $\l X_1^{\phi},V_1^{\phi} \r$, relying on theory of pseudo-differential operators. Fix $T>0$ and
consider the Banach space $L^1\l  \R^d\times [0,T] \r$ endowed with its natural norm. Consider 
\begin{align*}
    D := \lg g\in L^1\l  \R^d\times [0,T]\r  \mid g(x,0) = 0,\ \partial_t g \in L^1\l  \R^d\times [0,T]\r ,\ \nabla_x g \in L^1\l  \R^d\times [0,T]\r \rg.
\end{align*}
Let us call $\lg U_s^{cv},s\in\R_+ \rg$ the semigroup with generator $-\partial_t + cv\cdot \nabla_x$, that has domain containing $D$. The Phillips theorem (see \cite[Theorem 13.6]{schilling2012bernstein}) ensures that the subordinated semigroup $\{\mathcal{U}_s^{cv},\ s\in\R_+\}$ defined by the following Bochner integral  
\begin{align*}
    \mathcal{U}_s^{cv} g:= \int_0^\infty U_{\tau}^{cv} g P \l H(s) \in d\tau \r
\end{align*}
is again a strongly continuous contraction semigroup on $L^1\l  \R^d\times [0,T]\r$ and its generator is
\begin{align}
    -\phi(\partial_t - cv \cdot \nabla_x) g(y, s) &= \int_0^\infty (g(y + cvr, s - r)\mathds{1}_{\{r \le s\}} - g(y, s))\nu(dr), \label{Ricciuti toaldo: definizione operatore phillips}
\end{align}
defined at least on $D$, where $\nu$ is the Lévy measure associated to the subordinator $H$.

\begin{os} \label{Remark osservazione ricciuti toaldo operatore omogeneo}
    For $\phi(x) = x^\al$ the operator \eqref{Ricciuti toaldo: definizione operatore phillips} inherits the property of homogeneity from the $\al$-stable Lévy measure, if time and space are rescaled with the same parameter. Indeed, evaluating \eqref{Ricciuti toaldo: definizione operatore phillips}  in $(x/\gamma, t/\ga)$ and applying the change of variables $r = z/\gamma$ leads to 
    \begin{align*}
         -\l \l \partial_t - cv \cdot \nabla_x\r^\al g\r \l \frac{x}{\gamma},\frac{t}{\gamma} \r &= \int_0^\infty \l g\l\frac{x}{\gamma} + cvr, \frac{t}{\gamma} - r\r\mathds{1}_{\{r \le \frac{t}{\gamma}\}} - g\l \frac{x}{\gamma}, \frac{t}{\gamma}\r \r \frac{\al r^{-\al-1}}{\Gamma(1-\al)} dr \\
         &= \gamma^\al \int_0^\infty \l g\l\frac{x+cvz}{\gamma}, \frac{t-z}{\gamma} \r\mathds{1}_{\{z \le t\}} - g\l \frac{x}{\gamma}, \frac{t}{\gamma}\r \r \frac{\al z^{-\al-1}}{\Gamma(1-\al)} \gamma^{1+\al} \frac{dz}{\ga} \\
         &= -\ga^\al\l \partial_t - cv \cdot \nabla_x \r^\al \l g\l \frac{x}{\ga},\frac{t}{\ga} \r \r,
    \end{align*}
    for each $\ga \in \R_+$.
\end{os}

We are now ready for the following result, from which the equation for the density of $\l X_1^{\phi},V_1^{\phi} \r$ (if it exists) directly follows.

\begin{prop} \label{Proposizione equazione governante ricciuti toaldo}
    Let $x \mapsto h(x, v)$ be in $L^1(\mathbb{R}^d)$ for any $v$ and such that $\sup_v \|h(\cdot, v)\|_{L^1} < +\infty$ and $\sup_v \|v \cdot \nabla_x h(\cdot, v)\|_{L^1} < +\infty$. Denote
    \begin{align}
        q(x,v,t) := \mathbb{E}^v h\l X_1^\phi(t), V_1^\phi(t)\r. \label{Ricciuti toaldo evolution operator}
    \end{align}
    Then $q(x,v,t)$ solves the following equation 
    \begin{align}
        \phi(\partial_t - cv \cdot \nabla_x) q(x,v,t) =  \bar \nu (t)h(x+cvt, v) +  \lambda \int_{\mathcal{S}^{d-1}} (q(x,v',t) - q(x,v,t))\mu(dv'), \label{ricciuti toaldo equazione governante governing}
    \end{align}
    under the initial condition $q(x,v,0) = h(x,v)$,
    where $\mu$ denotes the uniform measure on $\mathcal{S}^{d-1}$ and the operator $\phi(\partial_t - cv \cdot \nabla_x)$ is defined in \eqref{Ricciuti toaldo: definizione operatore phillips}. Further, it satisfies the equation in the Fourier-Laplace space
    \begin{align}
        \phi(s + ick \cdot v) \widehat{\widetilde{q}}(k,v,s) = \frac{\phi(s + ick \cdot v)}{s + ick \cdot v} \widehat{h}(k,v) + \lambda \int_{\mathcal{S}^{d-1}} \l \widehat{\widetilde{q}}(k,v',s) - \widehat{\widetilde{q}}(k,v,s) \r\mu(dv'). \label{RT: fourier laplace equation for the density}
    \end{align}
\end{prop}

\begin{proof}
    For the reader convenience, we set $c = 1$. Up to time t, only two disjoint scenarios are possible: either there were no jumps or there was a jump at time $\tau\in (0,t)$. Using this, we can write the the renewal equation
    \begin{align}
        q(x,v,t) = \overline{\Psi}(t) h(x+vt, v) + \int_0^t \psi(\tau) \int_{\mathcal{S}^{d-1}} q(x+v\tau, v', t-\tau) \mu(dv') d\tau, \label{Equazione di rinnovo Ricciuti Toaldo}
    \end{align}
    where $\psi(\tau)$ represents the waiting time probability density function and $\overline \Psi(t) $ the survival function of the waiting time. 

    We apply the Fourier-Laplace transform $x\to k,\ t\to s$ to both sides of the renewal equation \eqref{Equazione di rinnovo Ricciuti Toaldo}. Let us evaluate the transform of each term individually.
    The free stream component has the following Laplace Fourier transform
    \begin{align*}
        \mathcal{L}\mathcal{F} \left\{ \overline{\Psi}(t) h(x+vt, v) \right\} &= \int_0^\infty e^{-st} \int_{\mathbb{R}^d} e^{ik \cdot x} \overline{\Psi}(t) h(x+vt, v) dx dt \\
        &=  \l  \int_0^\infty e^{-(s + ik \cdot v)t} \overline{\Psi}(t) dt  \r  \widehat{h}(k, v) = \widetilde{\overline{\Psi}}(s + ik \cdot v) \widehat{h}(k, v),
    \end{align*}
     after a change of variable $x' = x + vt$. 
     On the other hand, the scattering component has the following representation
     {\small
     \begin{align*}
         \mathcal{L}\mathcal{F} \left\{ \int_0^t \psi(\tau) \int_{\mathcal{S}^{d-1}} q(x+v\tau, v', t-\tau) \mu(dv') d\tau \right\} &= \int_0^\infty e^{-st} \int_{\mathbb{R}^d} e^{ik \cdot x}  \\& \qquad \l  \int_0^t \psi(\tau) \int_{\mathcal{S}^{d-1}} q(x+v\tau, v', t-\tau) \mu(dv') d\tau  \r  dx dt.
     \end{align*}}
     Exchanging the order of the integrals, the previous equation yields
     \begin{align*}
          & \int_0^\infty e^{-s \tau} e^{-ik \cdot v \tau} \psi(\tau)
        \int_{\mathcal{S}^{d-1}}  \l  \int_0^\infty e^{-st'} \int_{\mathbb{R}^d} e^{ik \cdot x'} q(x', v', t') dx' dt'  \r  \mu(dv') d\tau\\
         &=  \l  \int_0^\infty e^{-(s + ik \cdot v)\tau} \psi(\tau) d\tau  \r  \int_{\mathcal{S}^{d-1}} \widehat{\widetilde{q}}(k, v', s) \mu(dv')\\
         &= \widetilde{\psi}(s + ik \cdot v) \int_{\mathcal{S}^{d-1}} \widehat{\widetilde{q}}(k, v', s) \mu(dv').
     \end{align*}
     Substituting $\widetilde{\psi}(s) = \frac{\lambda}{\phi(s) + \lambda}$ and $\widetilde{\overline{\Psi}}(s) = \frac{1 - \widetilde{\psi}(s)}{s} = \frac{\phi(s)}{s(\phi(s) + \lambda)}$
    into the transformed renewal equation gives
     \begin{align*}
         \widehat{\widetilde{q}}(k, v, s) = \frac{\phi(s + ik \cdot v)}{(s + ik \cdot v)(\phi(s + ik \cdot v) + \lambda)} \widehat{h}(k, v) + \frac{\lambda}{\phi(s + ik \cdot v) + \lambda} \int_{\mathcal{S}^{d-1}} \widehat{\widetilde{q}}(k, v', s) \mu(dv')
     \end{align*}
     Multiplying both sides of the equation by $\phi(s + ik \cdot v) + \lambda$ and rearranging the equation one gets
     \begin{align*}
         \phi(s + ik \cdot v) \widehat{\widetilde{q}}(k, v, s) = \frac{\phi(s + ik \cdot v)}{s + ik \cdot v} \widehat{h}(k, v) + \lambda \int_{\mathcal{S}^{d-1}}  \l  \widehat{\widetilde{q}}(k, v', s) - \widehat{\widetilde{q}}(k, v, s)  \r  \mu(dv')
     \end{align*}
     so that Equation \eqref{RT: fourier laplace equation for the density} is proved. Furthermore, we observe that the Laplace-Fourier of the first lhs in \eqref{ricciuti toaldo equazione governante governing} is
     \begin{align*}
         &\int_0^{\infty} \int_{\R^d} e^{-st}e^{ik\cdot x} \phi\l \partial_t - v\cdot \nabla_x \r q(x,v,t) dx dt \\
         &= \int_0^{\infty} \int_{\R^d} e^{-st}e^{ik\cdot x} \l \int_0^{\infty}\l -q(x+v\tau , v, t-\tau) \mathds{1}_{\tau \leq t}  + q(x,v,t)\r \bar \nu (d\tau) \r dx dt \\
         &= \int_0^{\infty} e^{-st}\l \int_0^{\infty} \l -e^{-ik\cdot vt}\widehat{q}(k,v,t-\tau) \mathds{1}_{\tau \leq t} + \widehat{q}(k,v,t)\r \bar \nu (d\tau) \r dt 
     \end{align*}
    where we applied the change of variable $x' = x+vt$ which generates the term $e^{-ik\cdot vt}$ in the first addend and then integrated over $x'$. Separating the integrals, switching the order of integration in $t$ and $\tau$ and applying the change of variable $t'=t-\tau$, the previous equation yields
    \begin{align*}
        \widehat{\widetilde{q}}(k,v,s) \int_0^{\infty} \l 1 - e^{-\tau (s + i k \cdot v)} \r \bar \nu (d\tau)
         &= \widehat{\widetilde{q}}(k,v,s) \phi\l s + i k \cdot v \r
    \end{align*}
     and the thesis follows from uniqueness of Laplace-Fourier transform. 
\end{proof}

If $L$ is the inverse stable subordinator, the Equation \eqref{ricciuti toaldo equazione governante governing} reduces to \cite[Equation (5.16)]{ricciuti2023semi}.

\begin{os}
    We observe that, if the random flight starts in $(x,v) \in \R^d\times \S^{d-1}$, then $P\l X_1^\phi(t) \in \partial B_{ct}(x) \r \neq0$ implies that the distribution has a discrete component on the boundary of $B_{ct}(x)$ and there the density does not exist. Still, Equation \eqref{ricciuti toaldo equazione governante governing} for the evolution operator $q$ makes sense. Furthermore, whenever the density $f_t\l y,u \mid x,v \r := P \l X_1^\phi(t) \in dy, V_1^\phi(t) \in du \mid X_1^\phi(0) = x, V_1^\phi(0) = v \r$  exists in the open ball $\lg y\in \R^d \mid |y-x|_e < ct\rg$, it solves the forward equation 
    \begin{align}
        \phi \l \partial_t + cu \cdot \nabla_y\r f_t\l y,u \mid x,v \r =   \lambda \int_{\mathcal{S}^{d-1}} ( f_t\l y,u' \mid x,v \r - f_t\l y,u \mid x,v \r)\mu(du'), \label{FORWARD: ricciuti toaldo equazione governante governing}
    \end{align}
     Moreover, the following backward equation holds in $\lg y\in \R^d \mid |y-x|_e  \le ct\rg$
    \begin{align}
        \phi \l \partial_t - cv \cdot \nabla_x\r f_t\l y,u \mid x,v \r  = \bar \nu(t) \delta(y-x-cvt)\delta(u-v) +   \lambda \int_{\mathcal{S}^{d-1}} (f_t\l y,u \mid x,v' \r  - f_t\l y,u \mid x,v \r )\mu(dv'). \label{BACKWARD: ricciuti toaldo equazione governante governing}
    \end{align}
\end{os}

We now look for the governing equation of the law of the limiting process of $X_1^\phi$, under the hyperbolic scaling $(x,t)\mapsto(x/\eps,t/\eps)$. We restrict to $\phi(x)=x^\al$, $\al\in(0,1)$, so that we can exploit the property in Remark~\ref{Remark osservazione ricciuti toaldo operatore omogeneo} of  $\l\partial_t-cv\cdot\nabla_x\r^\al$. 

\begin{os}
    Let us consider the evolution operator $q(x,v,t)$ in \eqref{Ricciuti toaldo evolution operator}, an initial datum $g(y,v) = h(\eps y, v)$ for a small parameter $\eps>0$ and set $q^\eps(x,v,t) := q\l x/\eps, v, t/\eps \r$. Then 
    \begin{align*}
        q(x,v,t) = \E^v h\l \eps\l x+\int_0^t V_1^\phi(s)ds\r, V_1^\phi(t) \r \implies q^\eps(x,v,t) = \E^v g\l x+\eps\int_0^{t/\eps} V_1^\phi(s)ds, V_1^\phi\l \frac{t}{\eps} \r \r ,
    \end{align*}
    so that $q^\eps$ is the evolution operator of $\l x+\eps \l X_1^\phi(t/\eps) -x\r, V_1^\phi\l t/\eps\r \r$. In this sense, the scaling $(x,t)\mapsto(x/\eps,t/\eps)$ as $\eps \to 0$ coincides with $c\to \infty$ in \cite[(5.24)]{ricciuti2023semi}, hence we shall denote by $X_\infty(t)$ the limit (in distribution) of $ x+\eps\l  X_1^\phi\l t/\eps \r-x\r$. Such limit process is an anomalous diffusion with ballistic behaviour (for large times) and is constructed as sum between two terms: the first is a time changed rotationally invariant $\al$-stable process with the inverse of an $\al$-stable subordinator; the second is a corrective term, which accounts for the flight after the last velocity change. Details can be found in \cite[Theorem 5.9]{ricciuti2023semi}.
\end{os}

In \cite[Remark 5.7]{ricciuti2023semi} the governing equation of such limiting process is written for the distribution of the position at the last velocity change before time $t$; here we use Hilbert expansion to obtain the equation for the law for all times of the limit process, which also accounts for the flight in progress at time $t$. Let us denote by $\rho(x,t)$ the density of $X_\infty(t)$.

\begin{os} \label{Invarianza rotazioale si rho(x,t): osservazione ricciuti toaldo}
    The density $\rho$ is invariant under rotation, i.e.~$\rho(y,t) = \rho(-y,t)$ for each $y\in\R^d$. It directly follows from the symmetry of $\mu$. Indeed, the realization of velocities $\{ v_i, i\in\N\}$ which defines the trajectory $X_1^\phi(t)$ (hence $X_\infty(t)$ in the limit) has the same probability of $\{ -v_i, i\in\N\}$ which leads to $-X_\infty(t)$. 
\end{os}

\begin{te} \label{prop limiting fractional equation}
    Assume $q^\eps(x,v,t)=q_0(x,v,t)+\eps^\al q_1(x,v,t)+r(\eps)$,  where $q_0,q_1$ do not depend on $\eps$, for some remainder $r(\eps)$ such that, for each $(x,t) \in \R^d\times(0,\infty)$, $r(\eps)=o(\eps^\al)$ and $\l\partial_t-cv\cdot\nabla_x\r^\al r(\eps)=o(1)$ as $\eps\to0$. If the initial datum is independent of $v$, say $g(x,v) = g(x)$, then $q_0$ is independent of $v$. If also $g(x) = \delta_0(x)$, then $\rho(x,t) = q_0(x,t)$ and solves
\begin{align}
    \int_{\mathcal{S}^{d-1}} \l \l \partial_t + cv \cdot \nabla_x \r^\al \rho(x,t) - \frac{t^{-\al}}{\Gamma(1-\al)}\, \delta_0(x-cvt) \r \mu(dv) = 0. \label{Equazione governante diffusione anomala RT}
\end{align}
\end{te}

\begin{proof}
    We observe that $q^\eps(x,v,t) \to q_0(x,v,t)$ as $\eps \to 0$ and under the hypothesis on the initial datum 
    \begin{align*}
        q^\eps(x,v,t) = \E^v g\l x + \int_0^{t/\eps} V_1(s) ds \r \xrightarrow[]{\eps \to 0} \E g\l x+X_\infty(t) \r,
    \end{align*}
    so that $q_0(x,t) =  \E g\l x+X_\infty(t) \r$ which is independent of $v$. Moreover, by definition of expected value and using the invariance of $\rho(x,t)$ (see Remark \ref{Invarianza rotazioale si rho(x,t): osservazione ricciuti toaldo}) one gets
    \begin{align*}
        q_0(x,t) = \int_{\R^d} g(x+y) \rho(y,t)dy = \int_{\R^d} g(x-y) \rho(y,t) dy.
    \end{align*}
    If $g = \delta_0$, we obtain $q_0(x,t) = \rho(x,t)$ as required. We now aim to prove the governing equation for $\rho$. Evaluating Equation \eqref{ricciuti toaldo equazione governante governing} in $(x/\varepsilon, v, t/\varepsilon)$ and using Remark \ref{Remark osservazione ricciuti toaldo operatore omogeneo} one gets the rescaled equation
    \begin{align*}
        \eps^\al \l \l \partial_t - cv \cdot \nabla_x \r^\al q^\eps(x,v,t) -  \frac{t^{-\al}}{\Gamma(1-\al)} q\l \frac{x+cvt}{\eps},v,0 \r\r =    \lambda \l \mathcal{J}- \mathcal{I} \r q^\eps(x,v,t).
    \end{align*}
    Letting $\eps\to0$ gives $\l\mathcal{J}-\mathcal{I}\r q_0=0$, i.e.\ $\mathcal{J}q_0=q_0$ (according to what proved above, it does not depend on $v$).
    Instead, the terms of order $\eps^\al$ give
    \begin{align*}
        \l \partial_t - cv \cdot \nabla_x \r^\al q_0(x,t) - \frac{t^{-\al}}{\Gamma(1-\al)} \delta_0(x+cvt) = \lambda \l \mathcal{J}- \mathcal{I} \r p_1(x,v,t).
    \end{align*}
    Integrating over $v$ both sides, changing velocities $v\to -v$ (using that $\mu$ is symmetric) and using $q_0 = \rho$ yields  \eqref{Equazione governante diffusione anomala RT} as required.
\end{proof}
The use of Hilbert expansion method to get the anomalous diffusive equation for different Bernstein functions presents some difficulties, coming from the absence of property in Remark \ref{Remark osservazione ricciuti toaldo operatore omogeneo} of the non-local operators. In any case, I would argue that the attempt to extend this method represents a new line of research to be explored in future works.

\subsubsection{The telegraph equation} \label{Paragraph: ricciuti toaldo, telegrafo telegraph}
Here we present  the derivation of a time non-local counterpart of the telegraph equation \eqref{Introduzione equazione del telegrafo}.
Consider Equation \eqref{BACKWARD: ricciuti toaldo equazione governante governing} and  $p_+(x,t) := p(x,1,t), \ p_-(x,t) := p(x,-1,t)$. A substitution leads to the following system of coupled equations
\begin{align}
    \begin{cases}
        \phi \l \partial_t - c \partial_x  \r  p_+ = \bar \nu(t) \delta(y-x-ct) -\frac{\lambda}{2} (p_+ - p_-) \\ \phi \l \partial_t + c \partial_x \r  p_- = \bar \nu(t) \delta(y-x+ct)+ \frac{\lambda}{2} (p_+ - p_-) \label{ahshahshahsh}
    \end{cases}.
\end{align}
Now, consider the macro-observables $p(x,t) := p_+ + p_-$ and $w(x,t) := c(p_+ - p_-)$. Then the system \eqref{ahshahshahsh} can be written as
\begin{align*}
    \begin{cases}
        \phi \l \partial_t - c \partial_x  \r  \l p(x,t) + \frac{w(x,t)}{c} \r = 2\bar \nu(t) \delta(y-x-ct)   -\frac{\lambda}{c} w(x,t) \\
        \phi \l \partial_t + c \partial_x \r  \l p(x,t) - \frac{w(x,t)}{c} \r =  2\bar \nu(t) \delta(y-x+ct)  +  \frac{\lambda}{c} w(x,t)
    \end{cases}.
\end{align*}
Passing in Fourier-Laplace space $(x,t)\to (k,s)$ (recall that $\phi\l s + i c k   \r$ is the Laplace-Fourier symbol of $\phi \l \partial_t - c \partial_x  \r$) and summing and subtracting the equations one has 
{\small
\begin{align*}
    \begin{cases}
        \l \phi(s+ick) + \phi(s-ick) \r \widehat{\widetilde{p}}(k,s) +  \l \phi(s+ick) - \phi(s-ick) \r \frac{\widehat{\widetilde{w}}(k,s)}{c}  = 2e^{iky}\l \dfrac{\phi(s+ick)}{s+ick} + \dfrac{\phi(s-ick)}{s-ick}\r   \\
        \l \phi(s+ick) - \phi(s-ick) \r \widehat{\widetilde{p}}(k,s) +  \l \phi(s+ick) + \phi(s-ick) \r \frac{\widehat{\widetilde{w}}(k,s)}{c}  = -2 \frac{\lambda}{c} \widehat{\widetilde{w}}(k,s) + 2e^{iky}\l \dfrac{\phi(s+ick)}{s+ick} - \dfrac{\phi(s-ick)}{s-ick}\r 
    \end{cases}
\end{align*} }
Substituting $\widehat{\widetilde{w}}(k,s) /c$ from the first equation into the second and rearranging terms, one has 
\begin{align*}
&\lq\l \phi(s+ick)-\phi(s-ick)\r ^2 - \l \phi(s+ick)+\phi(s-ick)\r ^2 - 2\lambda\l \phi(s+ick)+\phi(s-ick)\r \rq \,\widehat{\widetilde p}(k,s) \\
&\qquad = -4e^{iky}\left[\phi(s+ick)\,\phi(s-ick)\l \frac{1}{s+ick}+\frac{1}{s-ick}\r  + \lambda\l \frac{\phi(s+ick)}{s+ick}+\frac{\phi(s-ick)}{s-ick}\r \right] \iff \\
&\phi(s+ick)\,\phi(s-ick)\,\widehat{\widetilde p}(k,s) + \frac{\lambda}{2}\l \phi(s+ick)+\phi(s-ick)\r \,\widehat{\widetilde p}(k,s) \\
&\qquad = e^{iky}\left[\phi(s+ick)\,\phi(s-ick)\l \frac{1}{s+ick}+\frac{1}{s-ick}\r  + \lambda\l \frac{\phi(s+ick)}{s+ick}+\frac{\phi(s-ick)}{s-ick}\r \right]
\end{align*}
Applying the inverse transform we obtain
\begin{align*}
    \phi(\partial_t + c\partial_x) \phi(\partial_t - c\partial_x) p(x,t) &= -\frac{\lambda}{2} \l \phi(\partial_t + c\partial_x) + \phi(\partial_t - c\partial_x) \r p(x,t) \\
    &+ \phi(\partial_t + c\partial_x)\lq \bar \nu (t)\,\delta(y-x-ct)\rq  + \phi(\partial_t - c\partial_x)\lq \bar \nu (t)\,\delta(y-x+ct)\rq  \\
    &+ \lambda\,\bar \nu (t)\l \delta(y-x-ct)+\delta(y-x+ct)\r,
\end{align*}
leading to the required semi-Markov telegraph equation for general Bernstein function $\phi$. In the case of $\phi(x) = x^\al,\ \al\in(0,1]$, it reduces to \cite[Equation (5.74)]{ricciuti2023semi} (with $\theta = \lambda/2$, $c = 1$ and considering that their formal derivation is on $|y-x|_e < t$).

\subsection{Distribution of the flight times and convergence to superdiffusion} \label{Paragraph: geometrically infinitely divisible distribution}

Let $W$ follow an exponential distribution. Then, if $J$ is one of the identically distributed flight times of $\l X_1^\phi,V_1^\phi \r$, one has $P(J>t) = P(W>L(t)) = P(H(W) > t)$ so that $J \overset{d}{=}H(W)$, being $H$ the subordinator, independent of $W$, with Bernstein function $\phi$. In this section we are interested in the following questions.
\begin{enumerate}
    \item Find characterization of the class of possible flight times, to which is devoted Subsection \ref{Subsection: H(W) coincide con la classe delle distribuzioni geometricamente infinitamente divisibili}.
    \item Among the class of possible flight times $H(W)$, find a characterization of $H$ such that $H(W)$ lies in the domain of attraction of a stable law, to which is devoted Subsection \ref{Subsection: ricciuti superdiffusion}. Within this class,  the limiting process of $X_1^\phi$ is a superdiffusion with ballistic behaviour, i.e.~the mean square displacement increases as $t^2$ for large $t$ (see \cite[Section 5]{ricciuti2023semi}).
\end{enumerate}
To the best of our knowledge, these characterizations are a novelty and are of independent interest in probability theory.

\subsubsection{On the distribution of the flight times} \label{Subsection: H(W) coincide con la classe delle distribuzioni geometricamente infinitamente divisibili}
For the reader convenience, we shall start by recalling some known facts.
A random variable $X$ is said to follow an infinitely divisible distribution iff $\forall n\in\N$ there exists a sequence $\lg X_i^{(n)},\ i\in\N_+ \rg$ of i.i.d.~random variables such that $X \overset{d}{=} \sum_{i=1}^n X_i^{(n)} $. Remarkable examples are the Normal, Gamma and Poisson distributions. The following result is standard: if $X$ takes values in $(0,\infty)$, then $X$ is infinitely divisible iff its Laplace transform is of the following form 
\begin{align*}
    \E e^{-sX} = e^{-h (s)}, \qquad s\geq 0, 
\end{align*}
where $h$ is a Bernstein function of the form \eqref{Rappresentazione Lévy Kintchine Bernstein function} with $a = 0$. The following subclass of infinitely divisible distributions has received considerable attention in probability theory (see \cite{klebanov1985problem} where it has been defined and \cite{klebanov1996sums} for a more general theory).
A random variable $X$ is said to be geometrically infinitely divisible if, for every $p \in (0, 1)$, there exists a sequence $\lg X_j^{(p)},\ j\in\N_+\rg$ of i.i.d.~random variables such that
\begin{align}
    X \overset{d}{=} \sum_{j=1}^{N_p} X_j^{(p)}, \label{Formula per la definizione di X geometricamente infinitamente divisibile}
\end{align}
where $P(N_p=k)=p(1-p)^{k-1}$, $k\in\N_+$, and $N_p$ is independent of
$\lg X_j^{(p)},\ j\in\N_+\rg$.
We have the following important result, which we present as corollary of \cite{fujita1993generalization}.

\begin{corollary} \label{Teorema: H(W) coincide con la classe delle distribuzioni geometricamente infinitamente divisibili}
    Let $\l X_1^\phi,V_1^\phi \r$ be the random flight defined in Section \ref{Ricciuti paragrafo}. Then the class of possible flight times coincides with the class of geometrically infinitely divisible distribution.
\end{corollary}

\begin{proof}
    The proof can be obtained by proving the following characterization of $X\in(0,\infty)$
    \begin{align}
        X \text{ is geometrically infinitely divisible } \iff \psi(s) :=\E e^{-sX} = \frac{1}{1+\phi(s)}, \quad s\geq0, \label{Formula: aratterizzazione di H(W) come infinit geometric divisibile}
    \end{align}
    for some Bernstein function $\phi(\cdot)$ of the form \eqref{Rappresentazione Lévy Kintchine Bernstein function} with $a = 0$ and $\lim_{x\to\infty}\phi(x) = \infty$. Such characterization is a consequence of \cite[Theorems 1.1, 1.4]{fujita1993generalization}.
    Indeed, we observe that its condition $\lim_{x\to0}f(x) = 0$ in \cite[(1.1)]{fujita1993generalization} implies $a = 0$ in \eqref{Rappresentazione Lévy Kintchine Bernstein function}. Moreover, if $b > 0$ then the condition $\lim_{x\to\infty} f(x) = \infty$ in \cite[(1.1)]{fujita1993generalization} directly follows. Instead, if $b = 0$ then $X \in (0,\infty)$ ensures that $\nu(0,\infty) = \infty$ and \cite[(1.1)]{fujita1993generalization} is again satisfied. 
    With this in mind, the proof follows from a Laplace transform argument. Indeed, 
    \begin{align*}
        \E e^{-s H(W)} = \E_W \lq \E \lq e^{-s H(W)}  \mid W = w\rq  \rq = \frac{1}{1 + \phi(s)}.
    \end{align*}
\end{proof}

\subsubsection{On the convergence to superdiffusion} \label{Subsection: ricciuti superdiffusion}
We recall that a positive monotone function $f$ on $(0,\infty)$ is said to be regular varying at infinity of index $\rho\in\R$ iff $f(\lambda x)/f(x) \xrightarrow[]{x\to \infty} \lambda^\rho$, $\forall \lambda \in \R_+.$
The class of regular varying functions of index $\rho$ is denoted by $\mathrm{RV}_\rho$. If $\rho = 0$, then we talk about slowly varying functions. Analogous definitions hold for regular variation at 0. 
Moreover, in this section we shall use the following definition. For an introduction to stable distributions, see e.g.~\cite[Section VI.1]{feller1991introduction}.
\begin{defin}
    Let $J$ be a random variable and let $\{J_k,\ k\ge 1\}$ be i.i.d.~copies of $J$. We say that $J$ (or its distribution) belongs to the domain of attraction of a non-degenerate random variable $R$ iff there exist sequences $\{a_n,\ n\ge 1\}\subset\R_+$ and $\{b_n,\ n\ge 1\}\subset\R$ such that
    \begin{align*}
        \frac{\sum_{k=1}^n J_k - b_n}{a_n} \xrightarrow[n\to \infty]{d} R,
    \end{align*}
    where $\xrightarrow{d}$ denotes convergence in distribution.
\end{defin}
We shall also speak of the domain of attraction of a distribution $F$, meaning the domain of attraction of a random variable with distribution $F$. A standard result is the following: if $\l \sum_{k=1}^n J_k - b_n \r/a_n$ converges in distribution to a non-degenerate limit, then the limit is stable. If $F$ is the law of $J$ and $\nu$ is a Lévy measure, we set $\overline{F}(x) := P(J > x)$ and $K_F(x) := \int_0^x y^2 F(dy),\ K_\nu(x) := \int_0^x y^2 \nu(dy)$. 

The goal of this section is prove the following theorem, which is one of the main results of this paper. 

\begin{te} \label{Teorema: dominio di attrazione stabile J = H(W)}
     Let $\l X_1^\phi,V_1^\phi \r$ be the random flight defined in Section \ref{Ricciuti paragrafo} and let $J$ be one of the i.i.d.~flight times, whose distribution is associated to a Bernstein function $\phi$ with Lévy measure $\nu$. Then for $\al \in (0,2)$, $J$ lies in the domain of attraction of an $\al$-stable law iff $\bar \nu (t) \in \mathrm{RV}_{-\alpha}$ as $t \to \infty$. In that case, $P(J > t) \sim \bar \nu (t)$. 
\end{te}

\begin{os} \label{Osservazione post Teorema dominio di attrazione: caso stabile spiegato nel dettaglio}
    Heuristically, the equivalence says that the randomization with the exponential time does not affect the heavy tails of $\nu$, meaning that $J$ is heavy-tailed exactly when $\nu$ is.  In the case of $\phi(s) = s^\al,\ \al\in(0,1)$, $J$ follows a Mittag-Leffler distribution, i.e.~$\E\lq e^{-sJ}\rq=1/(1+s^\al)$,  and scatterings occur at fractional Poisson time (see Remark \ref{Osservazione: subordinatore stabile e Caputo fractional derivative}). One has
    \begin{align*}
        \mathbb{E}\lq e^{-s \frac{\sum_{k=1}^n J_k}{n^{1/\alpha}}}\rq = \l \frac{1}{1 + \frac{s^\alpha}{n}}\r ^{n} = \l 1 + \frac{s^\alpha}{ n}\r ^{-n} \xrightarrow[]{n\to\infty} e^{-s^\alpha}, \qquad \al\in(0,1),
    \end{align*}
    which is the Laplace transform of an $\al$-stable distribution, hence the Mittag-Leffler distribution lies in the domain of attraction of an $\al$-stable law. Moreover, the Lévy measure of the $\al$-stable subordinator is $\nu(dy)=\frac{\al}{\Gamma(1-\al)}y^{-1-\al}dy$, so that 
    \begin{align*}
        \bar \nu(t) = \int_t^\infty \frac{\al}{\Gamma(1-\al)}y^{-1-\al}dy = \frac{t^{-\al}}{\Gamma(1-\al)}, \qquad t>0,
    \end{align*}
    which implies $\bar \nu \in \mathrm{RV}_{-\al}$. Finally, in this case  $P(J>t) = E_{\al} (-t^\al) \sim t^{-\al}/\Gamma(1-\al)$ for large times, where $E_\al(\cdot)$ is the Mittag-Leffler function, for $\al\in(0,1)$ (such asymptotic result is standard and can be obtained using Tauberian theorems \cite[Section XIII.5, Theorem 4]{feller1991introduction}).
    Theorem \ref{Teorema: dominio di attrazione stabile J = H(W)} extends this classical result to each subordinator whose Lévy measure has a regularly varying tail of index $-\al$, $\al\in(0,2)$.
\end{os}

Before discussing the proof of Theorem \ref{Teorema: dominio di attrazione stabile J = H(W)}, it is useful to present some other examples of subordinators $H$ satisfying the hypotheses and the corresponding distribution $J$.

\begin{ex}
    Assume that $H$ has Bernstein function $\phi(s) = \log \l 1+s^\al \r,\ \al\in(0,1)$ (see \cite[p. 328, No 59]{schilling2012bernstein}). In this case
    \begin{align*}
        \nu(dx) = \al x^{-1} E_{\al}\l -x^{\al} \r dx \sim \frac{\al}{\Gamma\l 1-\al \r} x^{-1-\al}dx \implies \bar \nu(t) \sim \frac{1}{\Gamma\l 1-\al\r} t^{-\al}, \quad t\to\infty,
    \end{align*}
    leading to $\bar \nu (t) \in \mathrm{RV}_{-\al}$ for $t\to\infty$. On the other hand, setting $J = H(W)$ as in Theorem \ref{Teorema: dominio di attrazione stabile J = H(W)}, 
    \begin{align*}
        \l \E e^{-\frac{s}{n^{1/\al}} J} \r^n = \lq 1+\log \l 1+\frac{s^\al}{n} \r \rq^{-n} \xrightarrow[]{n\to\infty} e^{-s^\al}, \quad \al\in(0,1),
    \end{align*}
    so that $J$ lies in the domain of attraction of an $\al$-stable law.
    In this case, the asymptotic behaviour coincides with the one obtained from the $\al$-stable subordinator as in Remark \ref{Osservazione post Teorema dominio di attrazione: caso stabile spiegato nel dettaglio}. However, $P(H(W) < \eps)$ differs as $\eps\to0$ so that they models different phenomena. 
\end{ex}

\begin{ex}
    Assume that $H$ has Bernstein function $\phi_1(s) = s\log \l 1 + s^{-1} \r $ (see \cite[p. 316, No 27]{schilling2012bernstein} with $a=1$). In this case
    \begin{align*}
        \nu_1(dx) = \frac{1-(1+x)e^{-x}}{x^2}\,dx, \quad x>0 \implies \bar \nu_1(t) =  \frac{1-e^{-t}}{t}\sim\frac{1}{t}, \quad t\to\infty,
    \end{align*}
    leading to $\bar \nu_1 \in \mathrm{RV}_1$. On the other hand, setting $J_k = H(W),\ k\in\N_+$ as in Theorem \ref{Teorema: dominio di attrazione stabile J = H(W)}, one has, after some manipulations,  
    \begin{align*}
        \E e^{-\frac{s}{n} \l \sum_{k=1}^n J_k - n\log n \r} = e^{s\log s + O\l n^{-1}\log^2 n \r} \xrightarrow[]{n\to\infty} e^{s\log s},
    \end{align*}
    which is the bilateral Laplace transform of an asymmetric $1$-stable distribution. Indeed, considering $Z$ and $Z_k,\ k\in\N$ with such distribution,
    \begin{align*}
        \mathbb E e^{-s(Z_1+\dots+Z_n)}=e^{ns\log s}=e^{ns\log s+ns\log n-sd_n} = \mathbb E e^{-s\l n^{1/\al}Z+d_n \r},
    \end{align*}
    for $d_n = n \log n$ and $\al=1$. Given that the set of Bernstein function is closed under composition (see \cite[Corollary 3.8]{schilling2012bernstein}), one may also consider $\phi_2(s) := \phi_1\l s^{\al} \r = s^\al \log\l 1+ s^{-\al} \r,\ \al\in(0,1)$. 
\end{ex}

It is useful to present also one case where Theorem \ref{Teorema: dominio di attrazione stabile J = H(W)} does not apply.

\begin{ex}
    Assume $H$ is the gamma subordinator, i.e.~it has Bernstein function $\phi(s) = \log(1+s)$ (see \cite[p. 314, No 26]{schilling2012bernstein} with $a=1$). In this case, 
    \begin{align*}
        \nu(dx) = x^{-1}e^{-x}dx, \quad x>0 \implies \bar \nu(t) \sim \frac{e^{-t}}{t}, \quad t\to\infty,
    \end{align*}
    so that the tail of the Lévy measure is not regularly varying. Setting $J = H(W)$, one has $\E J = 1$ and $\E J^2 = 3$ (using, e.g., the relation with Laplace transform) so that, if $J_k,\  k\in\N_+$ have the same distribution of $J$, 
    \begin{align*}
        \frac{\sum_{k=1}^n J_k - n} {\sqrt{2n}} \xrightarrow[]{n\to\infty} N(0,1),
    \end{align*}
    and $J$ lies in the domain of attraction of the normal distribution. 
\end{ex}

We now discuss some lemmas, from which the proof of Theorem \ref{Teorema: dominio di attrazione stabile J = H(W)} will follow. For $\phi(s)$ of the form \eqref{Rappresentazione Lévy Kintchine Bernstein function} with $a = 0$, we set 
$g(s) := -\phi''(s) = \int_0^\infty x^2 e^{-sx} \nu(dx)$. It shall be useful to observe that for $s\to 0$
\begin{align*}
    &\phi(s) \to 0, \quad \phi'(s) \to b+\int_0^\infty x\nu(dx), \quad g(s) \to \int_0^\infty x^2 \nu(dx),
\end{align*}
where we used dominated convergence theorem.
\begin{lem} \label{Ricciuti: Lemma AAAAA}
    Consider $J$ and $\nu$ as in Theorem \ref{Teorema: dominio di attrazione stabile J = H(W)}. Then 
    \begin{enumerate}[label=(\alph*)]
        \item If $\int_0^\infty x^2\,\nu(dx)<\infty$, then $ \mathbb E[J^2]<\infty$.
        \item If $\int_0^\infty x^2\,\nu(dx)=\infty$, then $\mathbb E\lq J^2e^{-sJ} \rq\sim g(s)$ as $s\to 0$.
    \end{enumerate}
\end{lem}

\begin{lem} \label{Ricciuti: Lemma BBBBB}
     Consider $J$ and $\nu$ as in Theorem \ref{Teorema: dominio di attrazione stabile J = H(W)}. Assume $\int_0^\infty x^2\nu(dx) = \infty$. Then $\forall \rho \geq0$
    \begin{align*}
        \int_0^x y^2 F(dy) \in \mathrm{RV}_{\rho} \iff \int_0^x y^2 \nu(dy) \in \mathrm{RV}_{\rho},
    \end{align*}
    and in that case $K_F(x)\sim K_\nu(x)$ as $x\to\infty$.
\end{lem}

It is worth recalling that, for every subordinator, $\int_0^1 x^2\,\nu(dx)\le\int_0^1 x\,\nu(dx)<\infty$, so that $\int_0^\infty x^2\,\nu(dx)=\infty$ if and only if $\int_1^\infty x^2\,\nu(dx)=\infty$; the finiteness of the second moment of $\nu$ depends only on its tail, i.e.~it is infinite precisely when large jumps of $H$ are frequent enough.
As the following example shows, it can be either finite or infinite so that the hypotheses of Lemmas \ref{Ricciuti: Lemma AAAAA} and \ref{Ricciuti: Lemma BBBBB} are non empty.   
\begin{ex}
    If $\phi(s) = \log(1+s)$, i.e.~$H$ is the gamma subordinator, then $\nu(dx) = x^{-1}e^{-x}dx$ and $\int_0^\infty x^2 x^{-1}e^{-x}\,dx < \infty$. Instead, if $\phi(s) = s^\al,\ \al\in(0,1)$, i.e.~$H$ is the $\al$-stable subordinator, one has $\nu(dx) = \frac{\al}{\Gamma(1-\al)} x^{-1-\alpha} dx$ and $\int_0^\infty x^2 x^{-1-\alpha} dx = \infty$.  
\end{ex}

\begin{lem} \label{Ricciuti: Lemma CCCCC}
    Let $m$ be a measure on $(0,\infty)$ with $\overline m(x):=m((x,\infty))<\infty$ and $K_m(x):=\int_0^xy^2\,m(dy)<\infty$ for all $x>0$, and let $\alpha\in(0,2)$. Then $\overline m\in \mathrm{RV}_{-\alpha}$ as $x\to\infty$ iff $K_m\in \mathrm{RV}_{2-\alpha}$ as $x\to\infty$, and in that case
    \begin{align}\label{Asintotica nel lemma strumentale per caratterizzazione nel bacino di attrazione}
        \lim_{x\to\infty}\frac{x^2\overline m(x)}{K_m(x)}=\frac{2-\alpha}{\alpha}.
    \end{align}
\end{lem}

Lemma \ref{Ricciuti: Lemma CCCCC} has a clear probabilistic interpretation. If the tail decreases as a power law of index $\al$, then the truncated second moment grows as $x^{2-\alpha}$ (giving infinite variance). In the limit $\al \to 2$, the contribution of the tails vanishes, recovering the normal domain. 

\begin{os}
    For the $\al$-stable Lévy measure, $\al\in(0,1)$, we already observed that $\nu(dy) =\frac{\al}{\Gamma(1-\al)} y^{-1-\al}dy$, which gives $\bar \nu(x) = \frac{x^{-\al}}{\al} $ and $\int_0^x y^2\nu(dy) = \frac{x^{2-\al}}{2-\al}$.
    As a consequence, the ratio $x^2\overline \nu /K_m(x)$ is constant and equal to $(2-\al)/\al$. This example is in a sense the guideline of Lemma \ref{Ricciuti: Lemma CCCCC}; indeed, regular variation heuristically coincides with an asymptotic self-similarity: what holds exactly for the tail of the $\al$-stable measure holds in the limit $x\to\infty$ for a regularly varying tail measure. 
\end{os}

Proofs of Lemmas \ref{Ricciuti: Lemma AAAAA},\ref{Ricciuti: Lemma BBBBB} can be found in Appendix \ref{Appendice Lemma AA, BB, CC}. Instead, Lemma \ref{Ricciuti: Lemma CCCCC} is \cite[Section VIII.9, Theorem 2]{feller1991introduction} with $\zeta = 2, \eta = 0$ (note that the hypothesis $U_2(\infty) = \infty$ is satisfied thanks to $\al\in(0,2)$).
We are now ready for the proof of the theorem.

\begin{proof}[Proof of Theorem \ref{Teorema: dominio di attrazione stabile J = H(W)}]
    From \cite[Section IX.8, Theorem 1a: alternative form]{feller1991introduction}, a non negative distribution $F$ belongs to the domain of attraction of an $\al$-stable law, $\al\in(0,2)$ iff 
\begin{align*}
        &1-F(x) \sim \frac{2-\al}{\al}x^{-\al} L(x), \quad \text{as $x\to \infty$} ,
\end{align*}
for some slowly varying function $L:(0,\infty) \to \R$.
Hence it suffices to prove that $\bar F(x) := P(J > x) \in \mathrm{RV}_{-\alpha} \iff \bar \nu(x)  \in \mathrm{RV}_{-\alpha}$ as $x\to\infty$. Assume $\bar \nu (x) \in \mathrm{RV}_{-\alpha}$. From Lemma \ref{Ricciuti: Lemma CCCCC} with $m = \nu$ (hypotheses are satisfied), one gets $K_\nu(x) \in \mathrm{RV}_{2-\alpha}$ and
        \begin{align*}
            \frac{x^2 \bar \nu (x)}{ K_\nu(x) }\to \frac{2-\al}{\al} \quad \text{as} \quad x\to \infty,
        \end{align*}
        which gives $K_\nu(x)\to\int_0^{\infty} y^2 \nu(dy) = \infty$. Lemma \ref{Ricciuti: Lemma BBBBB} ensures $K_F(x)\in \mathrm{RV}_{2-\alpha}$ and $K_F(x)\sim K_\nu(x)$ as $x\to \infty$. Using again Lemma \ref{Ricciuti: Lemma CCCCC} with $m=F|_{(0,\infty)}$ gives $\bar F(x)\in \mathrm{RV}_{-\alpha}$ and 
        \begin{align*}
            \frac{x^2 \bar F(x)}{ K_F(x) }\to \frac{2-\al}{\al} \quad \text{as} \quad x\to \infty.
        \end{align*}
        As a consequence,
        \begin{align*}
            \frac{\bar F(x)}{\bar\nu(x)}
            =\frac{x^2\bar F(x)/K_F(x)}{x^2\bar\nu(x)/K_\nu(x)}\frac{K_F(x)}{K_\nu(x)}
            \longrightarrow1\quad \text{as} \quad x\to\infty,
        \end{align*}
        as required. 

    Conversely, assume $\bar F(x)\in\mathrm{RV}_{-\al}$. Lemma \ref{Ricciuti: Lemma CCCCC} with $m=F|_{(0,\infty)}$ gives $K_F(x)\in\mathrm{RV}_{2-\al}$, hence $\E\lq J^2\rq=\lim_{x\to\infty}K_F(x)=\infty$, and point (a) of Lemma \ref{Ricciuti: Lemma AAAAA} yields $\int_0^\infty y^2\,\nu(dy)=\infty$. Then Lemma \ref{Ricciuti: Lemma BBBBB} gives $K_\nu (x) \in\mathrm{RV}_{2-\al}$ and $K_F(x)\sim K_\nu(x)$, and Lemma \ref{Ricciuti: Lemma CCCCC} with $m=\nu$ gives $\bar\nu (x)\in\mathrm{RV}_{-\al}$. Finally, $\bar F(x)\sim\bar\nu(x)$ follows exactly as above.

\end{proof}

The study of limiting (anomalous) diffusion for different subordinators, i.e.~those for which $H(W)$ is not in the domain of attraction of a stable law, is an open problem and can be object of future research.

\section{Non-Markovian transport processes: a time-change approach} \label{paragrafo de gregorio section}

In this section we introduce a second class of non-Markovian transport processes, showing anomalous diffusive behaviour and involving time non-local governing equations. To this end, consider the Markovian transport process $\l X^{(c,\lambda)}, V^{(c,\lambda)}\r$  recalled in \ref{Subsection: the markovian random flight}. Moreover, consider $L = \{L(t), t\geq0 \}$ the inverse of a subordinator $H$ characterized by the Bernstein function $\phi$ and Lévy measure $\nu$ and assume $H$ is independent of $\l X^{(c,\lambda)}, V^{(c,\lambda)}\r$. We use the notation $\l X_2^{\phi},V_2^{\phi} \r := \lg \l X_2^{\phi}(t), V_2^{\phi}(t)\r, \ t\geq0 \rg$ where 
\begin{align}
    X_2^{\phi}(t) := X^ {(c, \lambda)}(0) + c\int_0^{L(t)} V^ {(c, \lambda)}(s)ds, \qquad V_2^\phi(t) := \frac{dX_2^{\phi}(t)}{dt} \label{Random flight de gregorio definizione pos e velocità}
\end{align}
for the time-changed random flight. We do not make explicit the dependence on $c$ and $\lambda$ in order to simplify the notation. It is worth noting that $X_2^\phi$ does not evolve linearly in time so that such random motion does not match the Definition \ref{Definizione transport process random flight generica}.

\begin{os}
    From standard theory of time-change of Markov processes \cite{baeumer2009brownian,toaldo2015convolution} it is known that if a Markov process $\{Y(t),t\geq0\}$ has infinitesimal generator $G$, then the governing equation for $\E \lq  h(Y(L(t))) \mid Y(0) = x \rq$ is of the form $\phi\l \partial_t \r f(x,t) = G f(x,t)$ with $f(x,0) = h(x)$ (see also \cite{biovcic2026modelling} for a recent comprehensive survey). Such theory includes the case of $Y(t) = \l X(t), V(t) \r$ which would lead to the $\l X(L(t)), V(L(t)) \r$. Such motion differs from the one defined by \eqref{Random flight de gregorio definizione pos e velocità}. Indeed, $X_2^\phi(t) = X(L(t))$ but $\dot L(t) = 0$ so that $V_2^\phi(t) = 0$ almost surely. As we shall see, the mathematical treatment is similar but to keep a physical interpretation of the position-velocity  after the time-change we refer to \eqref{Random flight de gregorio definizione pos e velocità} instead of the process obtained by replacing the time to each coordinate with the inverse subordinator.
\end{os}

\begin{os}
    The random flight $\l X_1^{\phi},V_1^{\phi} \r$ models completely different phenomena compared to $\l X_2^{\phi},V_2^{\phi} \r$, when the particle remains trapped for random times.
    An example of relevant difference between the two models is the following: the random variable $X_1^\phi(t)$ has finite support. Indeed, from its explicit definition \eqref{ricciuti toaldo definizione random flight esplicita} one has, for $(x,v) \in \R^d \times \mathcal{S}^{d-1}$ initial position and velocity, $|X_1^\phi(t) - x| \leq ct$. In contrast, for each $M \in \R_+$ one has $P\l |X_2^\phi(t) - x |_e \geq M \r > 0$ for each $t\in\R_+$. Details can be found in \cite{angelani2024anomalous} in the case of $L$ the inverse stable subordinator, but the result holds in general. 
\end{os}

Section \ref{De Gregorio: governing equation equazione governante e proprietà distributizionali} is devoted to the discussion of the time non-local governing equation for the transport process and its anomalous diffusive behaviour, with explicit examples for different choice of $\phi$.
In Section \ref{paragraph paragrafo de gregorio anomalous diffusion} we deal with the hydrodynamic limit using Hilbert expansion techniques, which leads to the governing equation of the delayed  Brownian motion. In Section \ref{paragrafo: De Gregorio: limite lorentz} we prove that the non-Markovian transport process $\l X_2^\phi, V_2^\phi \r $ is an approximation of a Lorentz-type model, in a low-density regime.

\subsection{Governing equation and distributional properties} \label{De Gregorio: governing equation equazione governante e proprietà distributizionali}

Consider the evolution operator
\begin{align}
    q(x,v,t) := \mathbb{E}^{(x,v)} h\l X_2^\phi(t), V_2^\phi(t)\r = \mathbb{E}^{(x,v)} h\l X_2^\phi(t), 0\r, \label{operatore di evoluzione de gregorio}
\end{align}
for $h\in C_0 \l \R^d\times \mathcal{S}^{d-1}\r$. It is worth noting that, even though $V_2^\phi(t)$ does not depend on the initial velocity $v$, the random position $X_2^\phi(t)$ does so that we cannot ignore the dependence on $v$ of $q$. 
We further observe that such definition coincides with $q(x,v,t) = \E^{(x,v)} \Tilde{h} \l X(L(t)), V(L(t))\r$ under the choice $\Tilde{h} (x,v) = h(x,0)$. 
As a consequence, we present the governing equation for $q$ as a corollary of a more general theory, which has been recalled in Appendix \ref{Subsection: the markovian random flight}.

\begin{corollary}
    The map $t\to q(x,v,t)$ is the unique solution of the Cauchy problem
    \begin{align}
        \phi(\partial_t) u(x,v,t) = cv\cdot \nabla_x u(x,v,t) + \lambda \int_{\mathcal{S}^{d-1}} \l u(x,v',t) - u(x,v,t) \r \mu(dv'), \quad
            u(x,v,0) = h(x,0). \label{Equazione governante semigruppo De Gregorio}
    \end{align}
\end{corollary}
The formulation with the evolution operator is not only more general than the governing equation for the density of the random flight, but also the natural one, since the law of $\l X_2^\phi(t),V_2^\phi(t)\r$ does not even need to have a density.

We now aim to find the mean squared  displacement of such transport process, showing its anomalous diffusive behaviour. Assume $(x,v) \in \R^d \times \mathcal{S}^{d-1}$ is the starting point of the time-changed random flight.
We shall use $P^{x,\mu} := \int_{\S^{d-1}} P^{(x,v)} \mu(dv)$ and $\E^{x,\mu}$ the corresponding expectation. The mean squared displacement of the transport process is denoted by
\begin{align*}
    M_2  (t) := \mathds{E}^{x,\mu} \left | X_2^\phi(t) - x \right |_e^2, \qquad t\geq0.
\end{align*}
As for the random flight, we do not make explicit the dependence on $(c,\lambda)$.
\begin{prop} \label{Model 1: proposition on mean squared  displacement}
    Let $\l X_2^\phi,V_2^\phi \r$ be the non-Markovian transport process defined above. Then we have
    \begin{align*}
        M_2  (t)  = \frac{2c^2}{\lambda} S_{\lambda}' (z) \bigg|_{z=0} - \frac{2c^2}{\lambda^2} \l  1 - S_{\lambda}(t) \r ,
    \end{align*}
    where $S_{\lambda}$ is the survival function \eqref{funzione di sopravvivenza}.
\end{prop}
\begin{proof}
    The proof directly follows from a conditioning argument on $L(t)$ of $M_2(t)$ and the form of the mean square displacement for the Markovian random flights \eqref{Preliminaries: MSD del markovian random flight}. 
\end{proof}

Some considerations about the asymptotic behaviour are required in order to understand how $\l X_2^\phi,V_2^\phi \r$ behaves. We distinguish the cases $t\to 0^+$ and $t\to \infty$. 

We start with the analysis $t\to\infty$. It is true that $L(t) \xrightarrow[]{t\to\infty} \infty$, which gives $1-\E e^{-\lambda L(t)} \xrightarrow[]{t\to \infty} 0$. 
As a consequence, 
\begin{align*}
    M_2 (t) \sim \frac{2c^2}{\lambda} \E L(t) \qquad \text{as} \quad t\to \infty.
\end{align*}
\begin{ex}
    Assume $\phi(s)$ is regularly varying at $s\to 0$, so that $\phi(s) \sim s^{\alpha} l_0(1/s)$ as $s\to 0$, for some $\al \in [0,1]$ and $l_0(\cdot)$ slowly varying. Passing in the Laplace space $t\to s$, one has 
\begin{align*}
    \mathcal{L}\{\mathbb{E} L(t)\}(s) = \frac{1}{s\phi(s)} \implies \widetilde{M}_2(s) \sim  \frac{2c^2}{\lambda s^{\al + 1}l_0\l \frac{1}{s}\r} \qquad \text{as} \quad s\to 0.
\end{align*}
As a consequence, inverting $s\to t$ and using the Karamata Tauberian theorem \cite[Theorem XIII.5.4]{feller1991introduction} one gets
\begin{align*}
    M_2(t) \sim \frac{2c^2}{\lambda} \frac{1}{\Gamma (1+\alpha)} \frac{t^\al}{l_0(t)} \qquad \text{as} \quad t\to \infty.
\end{align*}
As a remarkable example, if $L$ is the inverse stable subordinator, i.e.~$\phi(s) = s^\al$ and $l_0(t) = 1$ for each $t\geq0$, one obtains 
\begin{align*}
    M_2(t) \sim \frac{2c^2}{\lambda \Gamma(1+\al)} t^\al \qquad \text{as} \quad t\to \infty,
\end{align*}
which is precisely the formula \cite[(3.9)]{angelani2024anomalous}. On the other hand, if one assumes $\phi(s) = \log(1+s)$, which corresponds to the gamma subordinator, the mean squared displacement grows linearly, i.e.
\begin{align*}
    M_2(t) \sim \frac{2c^2}{\lambda } t \qquad \text{as} \quad t\to \infty.
\end{align*}
It shows that, even if we time-change a diffusive process with the inverse of a subordinator, we may still get a diffusive behaviour as $t\to\infty$.
\end{ex}

Let us now investigate the case $t\to 0$. Recall that, referring to the  Lévy-Kintchine representation \eqref{Rappresentazione Lévy Kintchine Bernstein function}, we are considering subordinators that are not killed and strictly increasing.
Under these hypothesis, the  Karamata Tauberian theorem \cite[Theorem XIII.5.4]{feller1991introduction} applies and using
\begin{align}
        \widetilde{M}_2 (s) = \frac{2c^2}{s\phi(s) (\lambda + \phi(s))} = \frac{2c^2}{\lambda} \left[ \frac{1}{s\phi(s)} - \frac{1}{s\l \lambda + \phi(s)\r} \right]\label{MSD mean squared da invertire}
\end{align}
the following implication holds
\begin{align*}
    \widetilde{M}_2(s) \sim \frac{2c^2}{s\phi^2(s)}  \qquad \text{as} \quad s\to \infty \implies M_2(t) \sim \frac{2c^2}{\Gamma(1+2\al) \phi^2\l\frac{1}{t}\r} \to 0  \qquad \text{as} \quad t\to 0^+.
\end{align*}
How fast the mean squared displacement goes to zero as $t\to 0^+$ depends on the chosen subordinator. 
\begin{ex}
    If $\phi(s) = s^\al l_{\infty}(s)$ as $s\to \infty$ for some $\al\in [0,1]$ and slowly varying function $l_\infty(\cdot)$, then 
\begin{align*}
    M_2(t) \sim \frac{2c^2}{\Gamma(1+2\al)} \frac{t^{2\al}}{\l l_{\infty}\l\frac{1}{t}\r \r ^2}  \qquad \text{as} \quad t\to 0^+. 
\end{align*}
In the case of $\phi(s) = s^\al$, then 
\begin{align*}
    M_2(t) \sim \frac{2c^2}{\Gamma(1+2\al)} t^{2\al} \qquad \text{as} \quad t\to 0^+,
\end{align*}
recovering the behaviour \cite[(3.8)]{angelani2024anomalous}. Instead, for the gamma subordinator the mean squared displacement goes to zero slowlier, as the following implication shows
\begin{align*}
    l_{\infty}(s) = \log(s) \implies M_2(t) \sim \frac{2c^2}{\l \log\l\frac{1}{t}\r \r^2} \qquad \text{as} \quad t\to 0^+. 
\end{align*}
\end{ex}

\subsection{The limiting anomalous diffusion} \label{paragraph paragrafo de gregorio anomalous diffusion}

This section is devoted to different methods of deriving time non-local heat equations as a proper diffusion approximation of $\l X_2^\phi, V_2^\phi \r$. Among them, we shall discuss the Hilbert expansion method; as far as we know, it has not yet been applied to get time-fractional diffusion equations.

Firstly, we show that from standard arguments in probability theory, it can be proved that, under the scaling limit
\begin{align}
    c^2 \to \infty, \quad \lambda \to \infty \quad \text{s.t.} \quad \frac{c^2}{d\lambda} \to D \in \R_+, \label{de gregorio: hydrodynamic limit}
\end{align}
the random position $X_2^\phi$ converges weakly to the time changed Brownian motion. From a physical perspective, the average number of velocity changes increases, causing frequent collisions; at the same time, the particle speed grows to balance these rapid changes, leading to a non-trivial limit. Without loss of generality, we shall assume $D = 1$, $c^2 = d\lambda$ and let $c \to \infty$.
Let $A \in \R^d$ a Borel set with $\partial A$ (the boundary of $A$) of null Lebesgue measure. Conditioning on $L(t)$, an application of Corollary in \cite[page 169]{watanabe1970convergence} leads to 

\begin{align*}
    \lim_{c\to \infty} P \l X_2^\phi(t) \in A \r &= \int_0^\infty \lim_{c\to \infty} P \l X(u) \in A \r P(L(t) \in du) \\
    &= \int_0^\infty P \l B(u) \in A \r P(L(t) \in du) = P \l B(L(t)) \in A \r.
\end{align*}
Then from the Portmanteau lemma we get the convergence in distribution to the time-changed Brownian motion with the inverse subordinator $L$, being $L$ independent of the Brownian motion. As a consequence, we get the time non-local heat equation \eqref{Fractional heat equation - watanabe ed altri} as the governing equation for the law of the limiting (anomalous) diffusion. However, it is worth noting that Corollary in \cite[page 169]{watanabe1970convergence} cannot be applied in contexts as space non-homogeneous random flights, where $L(t)$ depends on the trajectory $\lg X(s),\ s\leq t \rg$.


Historically, one way to get the (anomalous) diffusion equation for the random position starting from the governing equation of a random flight, is to consider the case $d = 1$ which leads to the telegraph equation, and then apply an hydrodynamic limit \eqref{de gregorio: hydrodynamic limit}. We start by showing how this approach applies here.
In $d = 1$, the velocity vector simplifies to a scalar taking only $\pm 1$ 
and \ref{definizione operatore J di scattering} becomes $\mathcal{J}p(x,v,t) = \frac{1}{2} \left[ p(x,1,t) + p(x,-1,t) \right]$.

With the notation $p_+(x,t) := p(x,1,t), \ p_-(x,t) := p(x,-1,t)$, the forward version of Equation \eqref{Equazione governante semigruppo De Gregorio} gives us a coupled system of two equations
\begin{align} \label{de gregorio: Equazioni da sommare e sottrarre}
    \begin{cases}
        \phi \l \partial_t \r  p_+ = - c\partial_x p_+ - \frac{\lambda}{2} ( p_+ - p_-) \\ \phi \l \partial_t \r  p_- =  c\partial_x p_- + \frac{\lambda}{2} (p_+ - p_-)
    \end{cases}.
\end{align}
To decouple these equations and find the governing equation for the random position, we define the two macro observables $p(x,t) := p_+ + p_-$, which represents the total probability density, and $w(x,t) = c(p_+ - p_-)$ which plays the role of the flux. Summing the equations \eqref{de gregorio: Equazioni da sommare e sottrarre} yields the (time non-local) continuity equation
\begin{align}
    \phi \l \partial_t \r (p_+ + p_-) = -c\partial_x(p_+ - p_-) \iff \phi \l \partial_t \r p(x,t) = -\partial_x w(x,t),  \label{continuity equation de gregorio}
\end{align}
while subtracting the equations \eqref{de gregorio: Equazioni da sommare e sottrarre} and multiplying by $c$ yields the dynamic of the flux
\begin{align*}
    \phi \l \partial_t \r  c(p_+ - p_-) = -c^2\partial_x(p_+ + p_-) - \lambda c(p_+ - p_-) \iff \l \phi \l \partial_t \r + \lambda \r w(x,t) = -c^2 \partial_x p(x,t).
\end{align*}
To eliminate the flux $w(x,t)$, we apply the operator $ \l \phi \l \partial_t \r  + \lambda \r $ to \eqref{continuity equation de gregorio} to get
\begin{align*}
     \l \phi \l \partial_t \r  + \lambda \r  \phi \l \partial_t \r  p(x,t) = -\partial_x  \l \phi \l \partial_t \r  + \lambda \r  w(x,t).
\end{align*}
Substituting the expression for $ \l \phi \l \partial_t \r  + \lambda \r  w(x,t)$ into the right-hand side, we obtain
\begin{align}
\phi \l \partial_t \r ^2 p(x,t) + \lambda \phi \l \partial_t \r  p(x,t) = c^2 \partial_x^2 p(x,t). \label{De gregorio: fractional telegraph equation}
\end{align}
Equation \eqref{De gregorio: fractional telegraph equation} is the time non-local counterpart of the telegraph equation; the case $\phi(s) = s^\alpha,\ \alpha\in(0,1]$ reduces to \cite[Eq. 6.3]{angelani2024anomalous}. Under the limit \eqref{de gregorio: hydrodynamic limit}, the telegraph equation \eqref{De gregorio: fractional telegraph equation} (divided by $\lambda$) gives the required
\begin{align*}
    \phi(\partial_t) p(x,t) = \, \partial_x^2 p(x,t).
\end{align*}

We now discuss how the Hilbert expansion method applies here, leading to the time-fractional heat equation also in greater dimension and going beyond time-change techniques.

\subsubsection{A Hilbert expansion method leads to the time-fractional heat equation} \label{A Hilbert expansion method formally leads to the time-fractional heat equation}

We restrict the attention to $\phi(s) = s^{\alpha},\ \al\in(0,1]$, which corresponds to an $\al$-stable subordinator, whose properties were recalled in Remark \ref{Osservazione: subordinatore stabile e Caputo fractional derivative}. To this end, let us introduce a parameter $\varepsilon$ (small and positive) and consider the scaling 
\begin{align*}
    (x,t)\mapsto\l x/\eps,t/\eps^{2/\al} \r, \qquad q^\eps (x,v,t) := q\l \frac{x}{\eps}, v, \frac{t}{\eps^{2/\alpha}} \r,
\end{align*}
where $q$ is \eqref{operatore di evoluzione de gregorio} with initial datum $g(y,v) = h(\eps y, v)$. Note that it is a reformulation of the scaling \eqref{de gregorio: hydrodynamic limit}. Using homogeneity of the Lévy measure, such $q^\eps$ solves
\begin{align}
    \varepsilon^2  \partial_t ^{\al}  q^{\eps} (x,v,t) = \eps cv\cdot \nabla_x q^{\eps}  (x,v,t) +  \lambda  \l  \mathcal{J}-\mathcal{I}  \r  q^{\eps}  (x,v,t). \label{rescaled equation de gregorio}
\end{align} 
Consistently with the literature, we make the following hypothesis: $q^{\eps}(x,v,t)$ can be written as a sum with powers of $\eps$ as follows
\begin{align*}
    q^{\eps}(x,v,t) = q_0(x,v,t) + \eps q_1(x,v,t) + \eps^2 q_2 (x,v,t) + r(\eps),
\end{align*}
for some remainder $r(\varepsilon)$ that is $o\l\varepsilon^2\r$ as $\eps \to 0$.
Substituting the expression in the rescaled equation \eqref{rescaled equation de gregorio} we can group together the terms with the same power and compare terms of equal order.  The term with order $\eps^0$ gives 
\begin{align*}
    \lambda \l \mathcal{J} - \mathcal{I} \r q_0(x,v,t) = 0 \iff \J q_0(x,v,t) = q_0(x,v,t).
\end{align*}
On the other hand, we know that $\J p = \int_{\S^2}p(x,v',t)\mu(dv')$ does not depend on $v$  so that we set $q_0(x,v,t) =: \rho(x,t)$.
We now analyze the terms of order $1$. Comparing the coefficients, the obtained equation is 
\begin{align}
    -cv\cdot \nabla_x q_0(x,v,t) = \lambda \l \mathcal{J}-\mathcal{I} \r q_1(x,v,t) \implies  \lambda \l \mathcal{J}-\mathcal{I} \r q_1(x,v,t) = -cv\cdot \nabla_x \rho(x,t) .\label{Durante espansione di Hilbert per i termini di ordine eps} 
\end{align}
We aim to find the condition which $q_1(x,v,t)$ must satisfy.
Firstly, we observe that $\J - \Int$ is not invertible. Indeed, 
\begin{align}
    \ker \l\J - \Int\r = \lg f(x,v,t) \mid f(x,v,t) = h(x,t) \text{ for some } h(x,t) \rg \label{De gregorio limite diffusivo ker operatore}
\end{align}
hence the operator is not injective and its kernel is generated by the constant function $1$ in the space of velocities. For fixed $(x,t) \in \R^d \times (0,\infty)$, the operator $\mathcal{J}$ \eqref{definizione operatore J di scattering} acting on $v$ is linear and continuous by definition and is of finite-rank (the dimension of the image is 1), hence it is compact on $L^p\l\S^{d-1}\r$ (see \cite[Chapter 6, p. 157]{brezis2011functional}). 
Therefore, the Fredholm alternative \cite[Theorem 6.6, Remark 4]{brezis2011functional} applies with $T = \mathcal J$ and $f = \frac{c}{\lambda} v\cdot \nabla_x \rho(x,t)$: the equation \eqref{Durante espansione di Hilbert per i termini di ordine eps} has a solution iff the right-hand side $cv\cdot \nabla_x \rho(x,t)$ vanishes under the action of any functional in the kernel of the dual operator $(\mathcal{J}-\mathcal{I})'$, which acts on the dual space $L^q(\mathcal{S}^{d-1})$ (with $1/p + 1/q = 1$). 

The kernel of the dual operator is also generated by the constant function $1$ \cite[Theorem 6.6, point (d)]{brezis2011functional}. Testing the solvability condition yields
\begin{align*}
    \langle v \cdot \nabla_x \rho \, , \, 1 \rangle_{L^q} = \int_{\mathcal{S}^{d-1}} \l cv\cdot \nabla_x \rho(x,t) \r \mu(dv) = \l \int_{\mathcal{S}^{d-1}} v \, \mu(dv) \r \cdot \nabla_x \rho(x,t) = 0,
\end{align*}
because $\int_{\mathcal{S}^{d-1}} v \, \mu(dv)$ is the expected value of a symmetric (with respect to the origin) random variable. Consequently, the Fredholm alternative condition is satisfied, guaranteeing the existence of a solution. 

To get uniqueness, we look at the solution of the equation in the space of functions with null mean with respect to $v$. From a physical point of view, it is equivalent to saying that the total mass of the system is contained in $q_0(x,v,t)$. Indeed, $\rho(x,t) = \J q^{\eps}(x,v,t) \implies \J q_i(x,v,t) = 0,\ \forall i \in \N_+$. Under this assumption, Equation \eqref{Durante espansione di Hilbert per i termini di ordine eps} gives
\begin{align*}
    q_1(x,v,t) = \frac{c}{\lambda} v\cdot \nabla_x \rho(x,t).
\end{align*}
We now turn to the analysis of the second-order terms. The condition is 
\begin{align*}
    & \partial_t ^{\al} q_0(x,v,t) - cv\cdot \nabla_x q_1(x,v,t) = \lambda \l \mathcal{J} - \mathcal{I} \r q_2(x,v,t) \\
    &\implies 
    \partial_t ^{\al} \rho(x,t) -\frac{c}{\lambda} v\cdot \nabla_x \l cv\cdot \nabla_x \rho(x,t)\r = \lambda \l \mathcal{J} - \mathcal{I} \r q_2(x,v,t)
\end{align*}
Applying the operator $\mathcal{J}$ to each member and using that $\mathcal{J}$ is idempotent, one gets
\begin{align}
    \J \l \partial_t ^{\al} \rho(x,t) \r = \J \l \frac{c^2}{\lambda} v\cdot \nabla_x \l v\cdot \nabla_x \rho(x,t)\r \r \implies   \partial_t ^{\al} \rho(x,t)  = \J \l \frac{c^2}{\lambda} v\cdot \nabla_x \l v\cdot \nabla_x \rho(x,t)\r \r .\label{Equazione ad un certo punto: hilbert expansion}
\end{align}
We now observe what follows 
\begin{align*}
    v\cdot \nabla_x \l v\cdot \nabla_x \rho(x,t)\r = \sum_{i=1}^d \sum_{j=1}^d v_iv_j \frac{\partial^2 \rho(x,t)}{\partial x_i \partial x_j},
\end{align*}
where $v_i$ is the $i$-th component of the vector velocity $v\in \R^d$. Given that $v_iv_j$ is odd for $i\neq j$ and $\mu$ is the uniform measure on the sphere we have $J(v_iv_j) = 0$ for $i\neq j$. Instead, $J \l v_i^2 \r = 1 / d$ for each $i = 1,\ldots, d$. With this in mind the equation \eqref{Equazione ad un certo punto: hilbert expansion} becomes
\begin{align}
    \partial_t ^{\al} \rho(x,t)  = \Delta_x \rho(x,t). \label{de gregorio equazione di diffusione anomala anomalous heat equation}
\end{align}
Once more, this result is consistent with the literature.
When the small physical parameter $\eps$ is small but it cannot be assumed to be zero, the Hilbert expansion approach gives higher-order corrections to the anomalous heat equation  \eqref{de gregorio equazione di diffusione anomala anomalous heat equation}.

\subsection{The Lorentz model approximation} \label{paragrafo: De Gregorio: limite lorentz}

In this section we aim to prove that the non-Markovian transport process $\l X_2^\phi,V_2^\phi \r$ naturally arises as a scaling limit of a mechanical model of Lorentz type (see Appendix \ref{Preliminaries: gallavotti} for an introduction). We shall focus on dimension $ d = 3$; the case of $d = 2$ is the same. This section represents a central contribution of the paper, as it provides a physical interpretation for models and equations that are often introduced from a purely mathematical perspective.

In the Lorentz model recalled in \ref{Preliminaries: gallavotti},
the time between collision depends only on the geometry of the system.
Here, we extend this framework to more realistic scenarios, where the position of the obstacles does not suffice to explain the randomness of the motion.
In complex or porous media, the time between free flights is often prolonged by localized delays. 
To account for these random sojourn times, whose (possible) heavy-tailed distributions typically induce (subdiffusive) anomalous transport, we replace $t$ with a random operational time $L(t)$, defined as the inverse of a subordinator. This operational clock effectively pauses spatial propagation during trapping events, incorporating delay effects directly into the macroscopic limit.
A first consequence is the following: consider the cylinder of height $ct$ around the particle trajectory, capped by two hemispheres of radius $R$
\begin{align}
        \theta \l x,v,t\r := \left\{y\in \mathbb{R}^3 : | y-(x+cvs)|_e \leq R,  \,\,\,  s\in [0,t] \right\}.
    \end{align}
Then, under the low-density limit 
\begin{align}
    \rho \to \infty \quad   {\mathrm{and}} \quad R \to 0 \quad   {\mathrm{s.t.}} \quad  \rho c \pi R^2 \to \lambda \in (0,\infty) \label{low-density limit nel testo}
\end{align}
the free flight  time $T^R$ becomes 
\begin{align*}
    P\l T^R > t \r &= P \l N\l \theta (x,v,L(t) \r = 0 \r = \int_0^\infty g_t(d\tau) e^{-\rho \l \pi R^2 c \tau + \frac{4}{3}\pi R^3 \r} \\
    &\xrightarrow[]{R \to 0} \int_0^\infty g_t(d\tau) e^{-\lambda \tau} =  \E e^{-\lambda L(t)},
\end{align*}
which is precisely the survival function \eqref{funzione di sopravvivenza} for the waiting times of the underlying counting process. 

We are now ready to state the convergence theorem. We shall use the notation $P^{(x,v)} \l \cdot \r := P \l \cdot | \ \mathcal{X}^{\phi}(0) = x, \mathcal{V}^{\phi}(0) = v \r$.

\begin{te} \label{De gregorio: teorema gallavotti limite low-density bg}
    Let $L$ be the inverse process of a subordinator $H$ with Bernstein function $\phi$. Consider the position-velocity process  $\l \X^\phi , \V^\phi  \r :=  \lg \l \X^{\phi}(t), \V^{\phi}(t) \r,\ t\geq0 \rg$ in $\mathbb{R}^3 \times \mathcal{S}^2$, where $\X^{\phi}(t) := \X(L(t)) $ and $\V^{\phi}(t) :=  d\X(L(t))/dt$, with initial data given by the random variables $\l \X^\phi(0) , \V^\phi(0)  \r$ in $\mathbb{R}^3 \times \mathcal{S}^2$.
    For all times $t\geq0$, under the low-density limit \eqref{low-density limit nel testo},
    \begin{align*}
        \l \mathcal{X}^\phi(t),\mathcal{V}^\phi(t) \r \xrightarrow[R \to 0]{d} \l X_2^{\phi}(t), V_2^{\phi}(t)\r ,\
    \end{align*}
    where $\l X_2^{\phi}, V_2^{\phi} \r$ is the non-Markovian transport process introduced in Section \ref{paragrafo de gregorio section}. 
\end{te}
\begin{proof}
    For $h\in C_b(\mathbb{R}^3\times \mathcal{S}^2)$, let us write
\begin{align*}
    E^{(x,v)} h\l \X^\phi(t), \V^\phi(t) \r  = \int _{\R^3\times \mathcal{S}^2} h(x', v') P^{(x,v)} \l \mathcal{X}^{\phi}(t)\in dx', \mathcal{V}^{\phi}(t) \in dv' \r.
\end{align*}
To get the statement, we shall prove that
\begin{align}
    \lim _{R\to 0} E^{(x,v)} h\l \X^\phi(t), \V^\phi(t) \r  = E^{(x,v)} h \l X_2^\phi(t) , V_2^\phi(t)  \r \qquad \forall \, \,\,h\in C_b(\mathbb{R}^3\times \mathcal{S}^2) . \label{Model 2: Val atteso di h(pos ,vel)}
\end{align}
From a standard conditioning argument on $L(t)$ one has
\begin{align*}
    \lim _{R\to 0} E^{(x,v)}h\l \X^{\phi}(t), \V^{\phi}(t) \r  
    &= \lim_{R\to0} \int_0^{\infty} \mathds{E}^{(x,v)}h\l \X(\tau), 0 \r  P\l L(t) \in d\tau \r.
\end{align*}
Now we observe that $h$ is bounded and so is the expected value. Using the dominated convergence theorem, the previous equation yields
\begin{align*}
     \int_0^{\infty} \l \lim_{R\to0} \mathds{E}^{(x,v)}h\l \X(\tau), 0 \r \r P\l L(t) \in d\tau \r = \int_0^{\infty} \mathds{E}^{(x,v)}h\l X(\tau), 0\r P\l L(t) \in d\tau \r
\end{align*}
where we used the classical result by Gallavotti \cite{Gallavotti1972RigorousTO}.
The equality \eqref{Model 2: Val atteso di h(pos ,vel)} immediately follows.
\end{proof}

As in the case of \cite{Gallavotti1972RigorousTO}, the limiting procedure is such that the contribution of recollisions vanishes. However, as clear from Theorem \ref{De gregorio: teorema gallavotti limite low-density bg}, the limiting random flight is non-Markovian. Recent contributions discuss non-Markovian scaling limits for the Lorentz gas. Contrary to recent literature on Lorentz models \cite{bobylev2001liouville,facciaroni2025random,nota2024diffusion,nota2022two}, here non-Markovianity and non-local equations \eqref{Equazione governante semigruppo De Gregorio} arise from a time-change argument, while other characteristics of the system are the ones of \cite{Gallavotti1972RigorousTO}.

\section{Conclusion}

This paper aims to give a contribution in the study of non-Markovian random flights, focusing on two cases of interest, where the governing equation involve non-local operators.  

Going beyond the semi-Markov framework by considering some non-Markovian random flights which exhibit long-memory behaviour, due to stochastic dependence among their  waiting times, presents substantial technical challenges and likely explains the scarcity of contributions in the literature. However, they emerge naturally in physics it would be of interest to give a contribution in this direction. This section is devoted to possible open problems.

\subsection{Directions of research}

Consider the definition of random flight \ref{Definizione transport process random flight generica} under Assumption \ref{Assumption on random flights}. Then the random flight is uniquely characterized by the joint law of the flight times $(J_1,\dots,J_n)$, or, equivalently, by the underlying counting process $\mathcal{N}$. As a consequence, the random flight $(X,V)$ is non-Markovian exactly when $\mathcal{N}$ is, and in this case it inherits the properties of $\mathcal{N}$.

In \cite[Definition 2.3]{facciaroni2026non} the authors introduced a novel counting process, named exchangeable fractional Poisson process, which is fully non-Markovian and has infinite memory. It is characterized by the following survival function of the waiting times
\begin{align}
    P[J_1 > t_1,\ \ldots,\ J_n > t_n] = S\l  \sum_{k=1}^n t_k\r  \label{waiting times generalized fractional pp}
\end{align}
where $S(\cdot)$ is a completely monotone survival function, i.e.
\begin{align*}
    S(t) = \int_0^{\infty}e^{-st}\zeta(ds),
\end{align*}
for some probability distribution $\zeta$.
For $S(t) = E_\al(-t^\al)$ (which corresponds to $\zeta$ the Lamperti measure) the  counting process of \cite[Definition 2]{facciaroni2025markov} is re-obtained, where $E_\al(\cdot)$ is the Mittag-Leffler function. See
\cite{james2010lamperti} for details on the Lamperti distribution.

A novel (fully) non-Markovian random flight can be introduced by assuming fixed speed $c$ and the exchangeable fractional Poisson process \cite[Definition 2.3]{facciaroni2026non} as the counting process for the velocity-changes. It would be a natural generalization of the non-Markovian random flight introduced in \cite[Section 3]{facciaroni2025random}, where the authors assumed $S(t) = E_\al(-t^\al)$ (and thus \cite[Definition 2]{facciaroni2025markov} for the counting process). We conjecture that these random flights are the approximation, under a low-density limit, of Lorentz models where the particle collides elastically with spherical obstacles, whose centers are distributed according to a Poisson distribution of random intensity $\rho \mathcal{K}$, where $\mathcal{K} \sim \zeta$.

On the other hand, if the speed $\mathcal{C}$ is random, the random flight $(X,V)$ is identified by the joint law $\l \mathcal{C}, J_1,\ldots, J_n \r$ for each $n$. In \cite[Section 5]{facciaroni2025random} the authors introduced a random flight with random speed $\mathcal{C}$ Lamperti and the underlying counting process the exchangeable fractional Poisson process \cite[Definition 2]{facciaroni2025markov}, i.e.~with waiting times $J_k \overset{d}{=} \mathrm{Exp} \l \lambda /\Lambda \r$  where $\Lambda$ is Lamperti distributed. Such definition can be generalized by considering any other distribution $\zeta$ which is compatible with the definition of exchangeable fractional Poisson process \eqref{waiting times generalized fractional pp}. Such generalization is not trivial because if $\zeta$ is not Lamperti, then we cannot exploit the self-similar property anymore, which plays a fundamental role in \cite[Section 5]{facciaroni2025random}. We conjecture that the law of the resulting non-Markovian random flight solves the non-local equation
\begin{align*}
    \mathcal{D}_t^{\nu, -G} f(x,v,t) = - f(x,v,t), \qquad f(x,v,0) = h(x,t),
\end{align*}
where $G = -cv\cdot \nabla_x + \lambda \l\Int - \J \r$ and the operator $\mathcal{D}_t^{\nu, -G} $ is defined in \cite[(2.13)]{facciaroni2026non}.
As in the case of para-Markov chains \cite{facciaroni2026non}, we expect that the space and time operators do not de-couple, in general.

Finally, for each of the presented models would be useful to consider the non-isotropic case, even in non-euclidean spaces.

\appendix

\section{The Markovian random flight and the Lorentz model}

\subsection{Time-change of Markovian random flights} \label{Subsection: the markovian random flight}

Define the family of operators $\sgL^{\phi} := \lg \sgL_t^{\phi},\ t\geq0 \rg$ where 

\begin{align}
    \sgL_t^{\phi} h := \int_0^{\infty} \mathcal{T}^ {(c, \lambda)}_{\tau}\, h\, g_t(d\tau), \qquad h\in C_0\l \R^d \times \S^{d-1}\r, \label{operatore Q}
\end{align}
where $g_t$ is the measure of the inverse subordinator $L$. 
\begin{te}
    The set $\sgL^{\phi}$ is a family of uniformly bounded linear operators on $C_0\l \R^d \times \S^{d-1}\r$ and
    \begin{align*}
        \sgL^{\phi}_t : C_0\l \R^d \times \S^{d-1}\r \mapsto C_0\l \R^d \times \S^{d-1}\r.
    \end{align*}
    Furthermore, for any $h \in C_0\l \R^d \times \S^{d-1}\r$ the mapping $t \mapsto \sgL_t^{\phi}h$ is strongly continuous 
    and the map $t\to \sgL_t^{\phi} u_0$ is the unique solution of the Cauchy problem 
\begin{align}
    \phi(\partial_t) u(t) = \mathcal{A} u(t), \quad
            u(0) = u_0, \quad u \in D, \label{Problema di Cauchy astratto abstract problem}
\end{align}
    where $\phi\l\partial_t\r$ is defined in \eqref{operatore per semi markov} and $D = \lg h \in C_0 \l \R^d \times \S^{d-1} \r: c v \cdot \nabla h \in C_0 \l \R^d \times \S^{d-1} \r \rg$.
\end{te}
We re-write the proof of this known result for completeness. 
\begin{proof}
We divide the proof in three parts.
\begin{enumerate}
    \item The linearity of $\sgL^{\phi}$ simply follows from the definition \eqref{operatore Q}. Similarly, the fact that $\sgL^{\phi}_t h\in C_0\l \R^d \times \S^{d-1}\r$ for each $h\in C_0\l \R^d \times \S^{d-1}\r$ and $t\geq0$ comes from the definition of Bochner integral provided that the integral exists. This can be  seen by using \cite[Theorem 1.1.4]{arendt1987vector} that applies here since $\left\| \mathcal{T}(\tau)h  \right\| \leq \left\| h \right\|$. The same argument leads to contractivity, indeed,
       \begin{align}
        \left\|  \sgL^{\phi}_t h \right\| \, =\, & \left\| \int_0^{+\infty} \mathcal{T}^ {(c, \lambda)}_{\tau}\, h\, g_t(d\tau),  \right\| \notag \\
        \leq  \, & \int_0^{+\infty} \left\| \mathcal{T}^ {(c, \lambda)}_{\tau}\, h\,  \right\| \, g_t(d\tau) \notag \\
        \leq \, & \left\| h \right\|, 
    \end{align}
   which implies that $\lg \sgL^{\phi}_t \rg_{t\geq0}$ is uniformly bounded. 
   So, the operators $\lg \sgL^{\phi}_t \rg$ are also well-defined and contractive.

   \item We now prove strong continuity. Consider $s,t\geq0$. One has
   \begin{align*}
       \left\| \sgL_t^{\phi} h - \mathcal{Q}_s h \right\| &\leq \int_0^{\infty} \|\mathcal{T}^ {(c, \lambda)}_{\tau} h\| \, |g_t(\tau) - g_s(\tau)| \, d\tau \\
       &\leq  \|h\| \int_0^{\infty} |g_t(\tau) - g_s(\tau)| \, d\tau .
   \end{align*}
   As $t \to s$, $\int_0^{\infty} |g_t(\tau) - g_s(\tau)| \, d\tau \to 0$ (the integral is dominated by 2) and strong continuity follows. 
   \item We now show that the map $t\mapsto \sgL_t^{\phi} u_0= u(t)$ solves the Cauchy problem. Applying the Laplace transform $t\to s$ one has
   \begin{align*}
       \widetilde{u}(s) &= \int_0^{\infty} e^{-st} \l \int_0^{\infty} \mathcal{T}^ {(c, \lambda)}_{\tau}\, u_0\, g_t(\tau) d\tau \r dt \\
       &= \int_0^{\infty} \mathcal{T}^ {(c, \lambda)}_{\tau}\, u_0 \l \int_0^{\infty} e^{-st} g_t(\tau) dt \r d\tau \\
       &= \frac{\phi(s)}{s}\int_0^{\infty} e^{-\tau \phi(s) }\mathcal{T}^ {(c, \lambda)}_{\tau}\, u_0 d\tau,
   \end{align*}
   where we use the Fubini theorem for Bochner integrals (see \cite[Theorem 1.1.9]{arendt1987vector}) and the expression for the Laplace transform of the density of the inverse subordinator $\widetilde{g_s}(\tau) = - \partial_\tau \l e^{-\tau \phi(s)}/s\r = \phi(s) e^{-\tau\phi(s)}/s$ (the latter can be obtained using $P(L(t)\le\tau) = P(H(\tau) \ge t)$). The Laplace transform of the semigroup is exactly the resolvent (see e.g.~\cite[Theorem 1.9.1]{kolokoltsov2011markov} for generators and resolvents) so that one has 
   \begin{align*}
       \widetilde{u}(s) = \frac{\phi(s)}{s} \l \phi(s) I -\mathcal{A} \r^{-1} u_0.
   \end{align*}
   With this in mind, we apply the Laplace transform on both sides to the Equation \eqref{Problema di Cauchy astratto abstract problem}, which gives
   \begin{align*}
       \phi(s) \widetilde{u}(s) - \frac{\phi(s)}{s} u_0 =\mathcal{A} \widetilde{u}(s) \iff \widetilde{u}(s) = \frac{\phi(s)}{s} \l \phi(s) I -\mathcal{A} \r^{-1} u_0.
   \end{align*}
   The thesis follows from uniqueness of Laplace transform for functions taking values in Banach space.
\end{enumerate}
\end{proof}

\subsection{Low-density limit of Lorentz models} \label{Preliminaries: gallavotti}
The Markovian random flight may emerge as approximation of the motion of a particle among obstacles in random position \cite{Gallavotti1972RigorousTO,spohn1978lorentz}. We here recall the setting of \cite{Gallavotti1972RigorousTO}, where obstacles are identical hard spheres of radius $R$, whose centres are randomly distributed according to a Poisson point process of intensity $\rho \in \R_+$, meaning that 
\begin{align}
    P\l \bigcap_{j = 1}^n \lg \mathcal{N} (B_j)= k_j \rg\r = e^{-\sum_{j=1}^n \rho |B_j|}\prod _{j=1}^n  \frac{\l \rho |B_j |\r ^{k_j}}{k_j !} \label{Inhomogeneous Poisson field as a special case of fractional poisson field},
\end{align}
for each collection of mutually disjoint, finite  Borel sets $\{B_j,\ j \in\{1,\ldots, n\}\}$,  and any choice of non-negative integers $\{k_j,\ j \in\{1,\ldots, n\}\}$, being $\mathcal{N}(B)$ the random variable \textit{Number of points of the field inside $B$}. Some obstacles may overlap.
Outside the system of obstacles, the particle moves along a straight line with constant speed $c>0$. When reaching the surface of an obstacle, an elastic collision occurs. This means that the particle is specularly reflected and the post-collisional speed is again equal to $c$.
We denote by $ \mathcal{X}(t)$ the position of the particle at time $t$ and by $ \mathcal{V}(t)$ the unit velocity vector at time $t$, i.e.~the direction of the particle at time $t$.
The resulting position-velocity process  $\l \X , \V  \r := \lg \l \Xt, \Vt \r,\ t\geq0 \rg$ is the so called Lorentz process in $\mathbb{R}^3 \times \mathcal{S}^2$.
To simplify the notation, we will not make explicit the dependence on $\rho,R$.
A standard procedure to get ride of difficulties in studying the motion of the particle is to analyze the system in the low-density limit
\begin{align} \label{BG limit}
    \rho \to \infty \quad   {\mathrm{and}} \quad R \to 0 \quad   {\mathrm{s.t.}} \quad  \rho c \pi R^2 \to \lambda \in (0,\infty).
\end{align}
In this regime, the distribution of the free flight time converges to an exponential distribution with mean $\lambda ^{-1}$. Hence, condition  \eqref{BG limit} means that the limit procedure is such that  the free flight time of the particle maintains a strictly positive, finite mean. 
An important result of \cite{Gallavotti1972RigorousTO}, is that the Lorentz process $\l \X , \V  \r $ converges, under the low-density limit \eqref{BG limit}, to the Markovian random flight $\l X^{(c,\lambda)}, V^{(c,\lambda)} \r$; the main idea of the proof is that under the limit the contribution of recollisions vanishes leading to a Markov process (although the Lorentz process is fully non-Markovian). 

\

\section{Proofs of Lemmas \ref{Ricciuti: Lemma AAAAA},\ref{Ricciuti: Lemma BBBBB}} \label{Appendice Lemma AA, BB, CC}

\begin{proof}[Proof of Lemma \ref{Ricciuti: Lemma AAAAA}]
    We preliminary observe that, using the expression for the Laplace transform of $J$, one has
    \begin{align}
        \E \lq J^2 e^{-sJ} \rq &= \l \frac{1}{1+\phi(s)} \r'' = \frac{g(s)}{\l1+\phi(s)\r^2} + \frac{2 \l \phi'(s) \r^2}{\l 1+\phi(s) \r^3} \to \E \lq J^2 \rq. \label{Ricciuti: strumentale}
    \end{align}
    Let us divide the proof into two parts.
    \begin{enumerate}[label=(\alph*)]
        \item Using \eqref{Ricciuti: strumentale} yields $\E J^2 = g(0) + 2 \l b+\int_0^\infty x\nu(dx)\r^2 <\infty$. Indeed, for each subordinator $\int_0^1 x^2 \nu(dx) \leq \int_0^1 x\nu(dx) <\infty$ and, using the hypothesis, $\int_1^\infty x \nu(dx) \leq \int_1^\infty x^2 \nu(dx) < \infty$.
        \item Under $\int_0^\infty x^2\,\nu(dx)=\infty$ we have $g(s) \to \infty$ as $s\to 0$. From \eqref{Ricciuti: strumentale}, it suffices to show that $\l \phi'(s)\r^2 = o(g(s))$ as $s\to0$. Let us fix $M>0$. Then by Cauchy-Schwartz inequality (multiplying the integrand function by 1)
        \begin{align*}
            \int_M^\infty xe^{-sx}\nu(dx) \leq \l  \bar \nu(M) \int_M^\infty x^2 e^{-2sx}\nu(dx) \r^{\frac{1}{2}} \leq \sqrt{ \bar \nu(M) g(s) }.
        \end{align*}
        The second inequality holds because $x^2e^{-2sx}\le x^2e^{-sx}$.
        It follows that 
        \begin{align*}
            \phi'(s) &= b + \int_0^M xe^{-sx} \nu(dx) + \int_M^\infty xe^{-sx} \nu(dx) \leq b + \int_0^M xe^{-sx} \nu(dx) + \sqrt{ \bar \nu(M) g(s) }.
        \end{align*}
        As a consequence ($g(s)$ is positive), 
        \begin{align*}
            \frac{\l \phi'(s) \r^2}{g(s)} \leq \frac{ \l b + \int_0^M xe^{-sx} \nu(dx) \r^2}{g(s)} + \bar\nu(M) \implies \limsup_{s\to0} \frac{\l \phi'(s) \r^2}{g(s)}\le \bar\nu(M) .
        \end{align*}
        Letting $M\to 0$ one has $\l\phi'(s)\r^2 = o(g(s))$ for $s\to0$ and the thesis follows.
    \end{enumerate}
\end{proof}

\begin{proof}[Proof of Lemma \ref{Ricciuti: Lemma BBBBB}]
    We preliminarily observe that if $f(x)\in \mathrm{RV}_\rho$, then $x^{-\rho}f(x) \in \mathrm{RV}_0$. For this reason, if $K_F(x) \in \mathrm{RV}_\rho$, then $l(x) := \Gamma (1+\rho) x^{-\rho} K_F(x)$ is slowly varying. Recalling that $\int_0^\infty e^{-sx} K_F(dx) = \int_0^\infty e^{-sx}x^2 F(dx)$, we can apply \cite[XIII.5 Theorem 2]{feller1991introduction} with $U(x) = K_F(x),\ L = l$ and $\omega(s) = \int_0^\infty e^{-sx} K_F(dx) = \E \lq  J^2 e^{-sJ} \rq$, which gives
    \begin{align*}
        \E \lq J^2 e^{-sJ} \rq \sim s^{-\rho} l\l \frac{1}{s} \r, \quad s\to0 \implies g(s) \sim s^{-\rho} l\l \frac{1}{s} \r, \quad s\to0.
    \end{align*}
    The implication comes from \eqref{Ricciuti: strumentale}.
    Given that $g(s)$ is the Laplace transform of $K_\nu(x)$, using again \cite[XIII.5 Theorem 2]{feller1991introduction} with $U(x) = K_\nu(x)$ one has the thesis. The converse holds similarly.
\end{proof}

\section*{\small Acknowledgements}
The authors acknowledges financial support from Sapienza University of Rome, Fondo di Avvio alla Ricerca 2025 (Grant No. AR125199BF489F9E).

\section*{\small Conflict of interest and data statement} 

{\small The author declares that he has no conflicts of interest. The author declares that this is a theoretical paper and it does not make use of data.}

\printbibliography

\end{document}